\documentclass[a4j,10pt]{article}
\usepackage{graphicx}
\usepackage{color}
\usepackage{latexsym}
\usepackage{amsmath}
\usepackage{amsthm}
\usepackage{amssymb}
\newtheorem{thm}{Theorem}
\newtheorem{prop}[thm]{Proposition}
\newtheorem{claim}{Claim}

\newtheorem{lm}{Lemma}
\newtheorem{cor}[lm]{Corollary}

\newtheorem{conj}[thm]{Conjecture}

\newcounter{Case}
\newcounter{Subcase}[Case]
\def\komento#1{\ }

\def\qed{\hfill$\Box$}
\def\ed{d_H^e}
\def\F{\mathcal{F}}
\def\T{\Theta}

\def\Q{\mathcal{Q}}
\def\LC{{\mbox{\textit{LC}}}}
\def\LD{{\mbox{\textit{LD}}}}
\def\EL{{\mbox{\textit{EL}}}}
\def\cl{{\mbox{\textit{cl}}}}

\allowdisplaybreaks
\begin{document}
\hspace*{0.5cm}
\begin{center}
{\LARGE
Induced nets and Hamiltonicity of claw-free graphs}
\bigskip\\
{\large
Shuya Chiba\footnotemark[1] \footnotemark[2]\qquad
Jun Fujisawa\footnotemark[3] \footnotemark[4]\\
\footnotetext[1]{%
Department of Mathematics and Engineering,
Kumamoto University,
2-39-1, Kurokami, Kumamoto 860--8555,
Japan.
\texttt{schiba@kumamoto-u.ac.jp}
}%
\footnotetext[2]{%
work supported by Japan Society for the Promotion of Science,
Grant-in-Aid for Scientific Research (C) 17K05347}%
\footnotetext[3]{%
Faculty of Business and Commerce,
Keio University,
Hiyoshi 4--1--1, Kohoku-Ku,
Yokohama, Kanagawa 223--8521,
Japan.
\texttt{fujisawa@fbc.keio.ac.jp}
}%
\footnotetext[4]{%
work supported by Japan Society for the Promotion of Science,
Grant-in-Aid for Scientific Research (B) 16H03952 and
Grant-in-Aid for Scientific Research (C) 17K05349}%
}
\end{center}
\medskip
\begin{abstract}
The connected graph of degree sequence $3,3,3,1,1,1$ is called \textit{a net},
and the vertices of degree $1$ in a net is called its endvertices.
Broersma conjectured in 1993 that a $2$-connected graph $G$
with no induced $K_{1,3}$ is hamiltonian
if every endvertex of each induced net of $G$ has degree at least $(|V(G)|-2)/3$.
In this paper we prove this conjecture in the affirmative.
\end{abstract}
\vspace{0.5\baselineskip}
%
%
%
%
%
%
\section{Introduction}
\label{intro}
Hamiltonian cycles in graphs have been extensively studied in the literature
(cf., e.g., \cite{Go1, Go2, Go3}).
Several decades ago,
research was mainly focused on their relation to the four-color problem;
however, since the approval of Dirac's and Ore's theorems,
various studies have considered degree conditions.
While some of these degree conditions became extremely complicated
as this type of research progressed,
there remains several unsolved problems that can only
be stated in short, easily understandable form,
and these problems still engage our interest.
One of these is provided in Conjecture \ref{conj},
which generalizes two classical results obtained by Matthews and Sumner
(Theorem \ref{MS}) and by Duffus, Gould and Jacobson (Theorem \ref{DGJ}).
The connected graph of degree sequence $3,3,3,1,1,1$ is called \textit{a net},
and the vertices of degree $1$ in a net is called its endvertices.
Moreover, a graph is called \textit{claw-free} if it has no induced $K_{1,3}$.
\begin{conj}[Broersma \cite{conj}]
Let $G$ be a $2$-connected claw-free graph of order $n$.
If every endvertex of each induced net in $G$ has degree
at least $\frac{n-2}3$, then $G$ is hamiltonian.
\label{conj}
\end{conj}
\begin{thm}[\cite{MS}]
Let $G$ be a $2$-connected claw-free graph of order $n$.
If the minimum degree of $G$ is at least $\frac{n-2}3$, then $G$ is hamiltonian.
\label{MS}
\end{thm}
\begin{thm}[\cite{DGJ}]
Let $G$ be a $2$-connected claw-free graph.
If $G$ has no induced net, then $G$ is hamiltonian.
\label{DGJ}%
\end{thm}%
To verify that the degree condition in Conjecture \ref{conj} is sharp,
we assume that $B_1$, $B_2$ and $B_3$ are complete graphs of the same order
with $\{x_i, y_i, z_i\} \subseteq V(B_i)$.
Further, let $G$ be a graph obtained from $B_1 \cup B_2 \cup B_3$
by adding six edges $x_i x_j$, $y_i y_j$ ($1 \le i < j \le 3$).
We can observe that $G$ is non-hamiltonian,
and $\{x_i, z_i \mid 1 \le i \le 3\}$ induces a net
with endvertices $z_1, z_2$ and $z_3$ in $G$.
Since each $z_i$ has degree $|V(B_i)|-1 = \frac{|V(G)|-3}3$,
the degree condition is observed to be indeed sharp. 

The only partial solution to Conjecture \ref{conj}
that is known to the authors of this study
is a theorem by \v{C}ada et al.~\cite{CLNZ}, 
which states that Conjecture \ref{conj} is true if
the degree condition is strengthened to $\frac{n+5}3$.

In this article, we prove Conjecture \ref{conj}.
Our theorem relies heavily on the closure concept
that was introduced by Ryj\'a\v{c}ek \cite{Ry}.
In Section \ref{preliminaries}, we introduce the terminology %
and present the preliminary results related to Ryj\'a\v{c}ek's closure,
before we introduce some key lemmas in Section \ref{lemmas}.
In Section \ref{large}, we provide the proof of Conjecture \ref{conj}
for graphs that contain at least $33$ vertices.
Since the proof for the smaller graphs (at most $32$ vertices)
comprises a tedious case-by-case analysis,
which is not enlightening,
we have instead provided a sketch of the theorem in Section \ref{small}
with the complete theorem being provided in \cite{Arxiv}.
%
%
%
%
%
\section{Preliminaries}
\label{preliminaries}
For standard terminology and notation, we refer the readers to \cite{Ds}.
In this paper, a \textit{graph} or a \textit{simple graph} means a finite undirected graph
without loops or multiple edges.
A \textit{multigraph} may contain multiple edges but no loops.
For a graph $G$ and $v \in V(G)$,
$N_G(v)$ denotes the set of the neighbors of $v$,
and $d_G(v) = |N_G(v)|$.
Moreover, for $U \subseteq V(G)$,
$G[U]$ denotes the induced subgraph of $G$ induced by $U$.
For a graph $H$ and $x \in E(H)$,
let $N_H^e(x)= \{y \in E(H) \setminus \{x\} \mid
\mbox{$y$ is adjacent to $x$ in $H$} \}$
and $\ed(x) = |N_H^e(x)|$.
The set of the endvertices of $x$
is denoted by $V_H(x)$, or simply $V(x)$.
We denote the set of vertices of degree one in $H$ by $V_1(H)$,
and a \textit{pendant edge} is an edge in which one endvertex has degree one.
For a vertex $v \in V(H)$,
the set of all the pendant edges which are incident with $v$
is denoted by $l_H(v)$, or simply $l(v)$.
For $X \subseteq V(H)$ and $e \in E(H)$,
we say that $e$ is dominated by $X$
if $V(e) \cap X \neq \emptyset$.
We often identify a subgraph $H$ of a graph $G$
with its vertex set $V(H)$.
The complete bipartite graph $K_{1,3}$ is called \textit{a claw},
and \textit{a clique} is a maximal complete subgraph of a graph.

In the rest of this section,
we prepare previous studies which are commonly used
in hamiltonian graph theory for claw-free graphs.
Let $G$ be a claw-free graph.
We call a vertex $v$ of $G$ \textit{locally connected}
(resp.~\textit{locally disconnected})
if $G[N_G(v)]$ is connected (resp. disconnected).
For a locally connected vertex $v$ of $G$,
the operation of joining all pairs of nonadjacent vertices in $N_G(v)$
is called \textit{the local completion at $v$}.
In \cite{Ry}, it is shown that this operation preserves
the claw-freeness of the original graph.
Iterating local completions,
we obtain a graph $G^*$ in which $G^*[N_{G^*}(v)]$ is a complete graph
for every locally connected vertex $v$.
We call this graph the \textit{closure} of $G$,
and denote it $\cl(G)$.
The closure of a claw-free graph has the following properties.
\begin{thm}[Ryj\'a\v{c}ek \cite{Ry}]
\label{Ry}
Let $G$ be a claw-free graph.
Then
$\cl(G)$ is uniquely defined and is the line graph of some triangle-free simple graph.
Moreover, $G$ is hamiltonian if and only if $\cl(G)$ is hamiltonian.
\end{thm}

For a claw-free graph $G$,
the set of locally connected vertices
and locally disconnected vertices
are denoted by $\LC(G)$ and $\LD(G)$,
respectively.
If $v \in \LC(G)$ and $G[N_G(v)]$ is not a complete graph,
then we call $v$ \textit{an eligible vertex},
and let $\EL(G)$ denote the set of eligible vertices of $G$.
For $x \in \EL(G)$,
$G_x$ denotes the graph obtained from $G$
by local completion at $x$.
It is shown in \cite[Lemma 9]{Pf} that
$\LD(G_x) \subseteq \LD(G)$ holds
for every claw-free graph $G$ and $x \in \EL(G)$.
This yields the following.
\begin{prop}
Let $G$ be a claw-free graph and let $v \in \LC(G)$,
then $v \in \LC(\cl(G))$.%
\label{lc}%
\end{prop}%

The following theorem is a basic tool
for the study on the hamiltonicity of line graphs.
A closed trail $T$ in a graph $H$ is called \textit{a dominating closed trail},
or a DCT,
if every edge of $H$ is dominated by $T$.
\begin{thm}[Harary and Nash-Williams \cite{HN}]
Let $H$ be a multigraph with $|E(H)| \ge 3$.
Then the line graph $L(H)$ is hamiltonian
if and only if $H$ has a DCT.
\label{HN}
\end{thm}
Let $H$ be a multigraph and $F$ be a subgraph of $H$.
Then $H/F$ denotes the multigraph
obtained from $H$ by identifying the vertices of $F$
with a new vertex, which is denoted by $v_F$,
and by deleting the created loops.
A multigraph $H$ with at least two vertices
is called \textit{collapsible}
if for every $S \subseteq V(H)$ with $|S|$ even,
there exists a spanning connected subgraph $F$ of $H$ such that
$v \in S$ if and only if $d_F(v)$ is odd.
Collapsible subgraphs have the following property.
%
\begin{prop}[Catlin \cite{Ca}]
Let $H$ be a multigraph and $F \subset H$ be a collapsible subgraph.
\begin{itemize}
\item[\normalfont{i)}]
If $H/F$ has a DCT containing $v_F$,
then $H$ has a DCT containing all the vertices in $F$.
\item[\normalfont{ii)}]
If $H/F$ is collapsible, then $H$ is collapsible.
\end{itemize}
\label{col}
\end{prop}

By $\F$ we denote the class of graphs
obtained by taking two vertex-disjoint triangles
$a_1 a_2 a_3$ and $b_1 b_2 b_3$
and by joining every pair of vertices $\{a_i, b_i\}$
by a path of length at least two
or by a triangle.
\begin{thm}[Brousek \cite{Br}]
Every non-hamiltonian $2$-connected claw-free graph
contains an induced subgraph
which is in $\F$.
\label{Br}
\end{thm}
%
%
%
%
%
%
%
\section{Lemmas}
\label{lemmas}
In this section we prove some lemmas.
Among them, Corollary \ref{endmove} (which is derived from Lemma \ref{em_g})
plays an important role in our proof.
We denote an induced net of $G$ with six vertices $x_1, x_2, x_3, y_1, y_2, y_3$
by $N(x_1, x_2, x_3; y_1, y_2, y_3)$
if $\{x_i x_{i+1}, x_i y_i \mid 1 \le i \le 3\} \subseteq E(G)$,
where $x_4 = x_1$.
\begin{lm}
Let $G$ be a claw-free graph and
let $G^0, G^1, \ldots , G^{l-1}, G^l$ be a sequence of graphs
such that $G^0 = G$, $G^l = \cl(G)$ and
$G^i$ is obtained from $G^{i-1}$ by the local completion of
an eligible vertex of $G^{i-1}$, for $1 \le i \le l$.
Moreover, let $N(x_1, x_2, x_3; y_1, y_2, y_3)$ be an induced net of $\cl(G)$,
let $R_0$ be the clique of $\cl(G)$ which contains the triangle $x_1 x_2 x_3$,
and let $R_j$ be the clique of $\cl(G)$ which contains the edge $x_j y_j$
for $1 \le j \le 3$.
Then the following holds:
\begin{itemize}
\item[\normalfont{i)}]
For each $i$ with $0 \le i \le l$,
there exists an induced net
$N^i = N(x_1^i, x_2^i, x_3^i; y_1^i, y_2^i, y_3^i)$ of $G^i$ such that
$x_j^i \in \{x_j\} \cup (V(R_0) \cap \LC(\cl(G)))$ and
$y_j^i \in \{x_j, y_j\} \cup (V(R_0 \cup R_j) \cap \LC(\cl(G)))$ for $j=1,2,3$.
\item[\normalfont{ii)}]
$x_1 x_2 x_3$ is a triangle in $G$ if and only if
$y_j^0 \in \{y_j\} \cup (V(R_j) \cap \LC(\cl(G))$ for each $j$.
\end{itemize}
\label{em_g}
\end{lm}
\noindent
\textit{Proof.}\quad
First we prove i) by reverse induction on $i$.
If $i=l$, then $N(x_1, x_2, x_3; y_1, y_2, y_3)$ is the desired induced net.
Assume that
$N^i = N(x_1^i, x_2^i, x_3^i; y_1^i, y_2^i, y_3^i)$
is the desired induced net of $G^i$ for some $i$ with $1 \le i \le l$.
Let $F_i = E(N^i) \setminus E(G^{i-1})$.
If $F_i = \emptyset$,
then $N^i$ is also the desired induced net of $G^{i-1}$,
and hence we assume that
$F_i \neq \emptyset$.
Let $u \in \EL(G^{i-1})$ such that
$G_u^{i-1} = G^i$.
Then,
all the edges of $F_i$ are contained in the clique induced by $N(u)$ in $G^i$.
Hence either
$F_i = \{x_j^i y_j^i\}$ for some $j$
or
$F_i \subseteq \{x_1^i x_2^i, x_2^i x_3^i, x_3^i x_1^i\}$.

First we consider the former case.
We may assume without loss of generality that $j=1$.
Let $y_1^{i-1} =u$, $y_j^{i-1} = y_j^i$
for $j = 2,3$
and $x_j^{i-1} = x_j^i$ for $1 \le j \le 3$.
Since $x_1^i y_1^i \in F_i$, we have $ux_1^i \in E(G^{i-1})$.
Moreover, since the neighbors of $u$ in $G^{i-1}$ induce a complete graph in $G^i$,
$u v \notin E(G^{i-1})$ for each $v \in \{x_2^i, y_2^i, x_3^i, y_3^i\}$.
Hence
$\{x_1^{i-1}, x_2^{i-1}, x_3^{i-1}, y_1^{i-1}, y_2^{i-1}, y_3^{i-1}\}$
induces a net, say $N^{i-1}$, in $G^{i-1}$.
Recall that $N^i$ satisfies i).
Since $u x_1^i y_1^i$ is a triangle in $G^i$,
it is a triangle in $\cl(G)$ as well,
and hence 
$u$ and $y_1^i$ are contained in
the same clique ($R_0$ or $R_1$) of $\cl(G)$.
Moreover, by Proposition \ref{lc},
$u \in \LC(G^{i-1}) \subseteq \LC(\cl(G))$.
Thus $N^{i-1}$ is the desired induced net of $G^{i-1}$.

Next we consider the case
$F_i \subseteq \{x_1^i x_2^i, x_2^i x_3^i, x_3^i x_1^i\}$.
If $|F_i| = 3$, then $\{u, x_1^i, x_2^i, x_3^i\}$
induces a claw in $G^{i-1}$, a contradiction.
Moreover, if $|F_i| = 1$,
say $F_i = \{x_1^i x_2^i\}$,
then $\{x_3^i, x_1^i, x_2^i, y_3^i\}$ induces a claw in $G^{i-1}$,
a contradiction.
Therefore $|F_i| = 2$, say $F_i = \{x_1^i x_2^i, x_2^i x_3^i \}$.
Then, since the neighbors of $u$ in $G^{i-1}$ induce a complete graph in $G^i$,
we have $u x_j^i \in E(G^{i-1})$ and
$u y_j^i \notin E(G^{i-1})$ for $j=1,2,3$.
Moreover, we have $x_2^i y_j^i \notin E(G^{i-1})$ for $j=1,3$,
since $N^i$ is an induced net in $G^i$.
Let $x_2^{i-1} = u$ and $y_2^{i-1} = x_2^i$.
Furthermore, for $j=1,3$, let
$x_j^{i-1} = x_j^i$ and $y_j^{i-1} = y_j^i$.
Then
$\{x_1^{i-1}, x_2^{i-1}, x_3^{i-1}, y_1^{i-1}, y_2^{i-1}, y_3^{i-1}\}$
induces a net, say $N^{i-1}$, in $G^{i-1}$.
Recall that $N^i$ satisfies i).
Since $u x_1^{i-1} x_3^{i-1}$ is a triangle in $G^{i-1}$,
it is a triangle in $\cl(G)$ as well,
and hence $u$, $x_1^{i-1}$ and $x_3^{i-1}$ are contained in $R_0$.
Moreover, by Proposition \ref{lc},
$u \in \LC(G^{i-1}) \subseteq \LC(\cl(G))$,
and by the induction hypothesis,
$y_2^{i-1} = x_2^i \in \{x_2\} \cup (V(R_0) \cap \LC(\cl(G)))$.
Thus $x_2^{i-1}, y_2^{i-1}$ satisfies i), and hence
$N^{i-1}$ is a desired induced net of $G^{i-1}$.
Consequently, i) holds for each $i$ with $0 \le i \le l$.

By the above procedure,
$x_1 x_2 x_3$ is a triangle in $G$
if and only if $|F_i| \le 1$ for every $1 \le i \le l$,
and $y_j^0 \in \{y_j\} \cup (V(R_j) \cap \LC(\cl(G))$ for each $j$
if and only if $|F_i| \le 1$ for every $1 \le i \le l$.
Hence ii) holds.
\qed
\medskip

The graph with seven vertices which is obtained from a claw
by subdividing each edge once is called \textit{a subdivided claw},
and the vertex of degree three in a subdivided claw is called
its \textit{center}.
Moreover, we call an edge $y$ of a graph $H$ \textit{heavy}
if $\ed(y) \ge \frac{|E(H)|-2}3$.

\begin{figure}[tb]
\begin{center}
\begin{picture}(120,100)
\put(10,22){\footnotesize$R_0$}
\put(50,55){\footnotesize$R_1$}
\put(50,22){\footnotesize$R_2$}
\put(50,-8){\footnotesize$R_3$}
\put(-34,63){\footnotesize$V_1(H) \ni$}
\put(41,96){\footnotesize$\in V_1(H)$}
\put(25,52){\footnotesize$x_1$}
\put(62,70){\footnotesize$y_1$}
\put(80,55){\footnotesize$R_1^+$}
\put(80,22){\footnotesize$R_2$}
\put(80,-8){\footnotesize$R_3$}
\includegraphics[width=3cm]{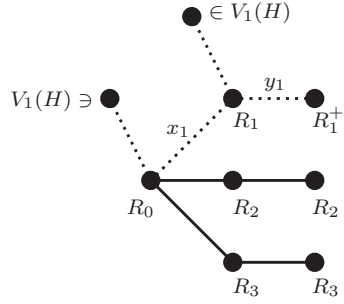}
\end{picture}
\caption{Candidates for $y_1'$ (dotted edges)}
\label{fig_endmove}
\end{center}
\end{figure}

In the rest of this section and the next section,
when we consider a claw-free graph $G$
and a triangle-free graph $H$ such that $L(H) = \cl(G)$,
each vertex of $H$ is denoted by a capital letter.
\begin{cor}
Let $G$ be a $2$-connected claw-free graph of order at least $3$
such that
every endvertex of each induced net in $G$ has degree
at least $\frac{|V(G)|-2}3$.
Let $H$ be the triangle-free graph such that $L(H) = \cl(G)$,
and let $\Lambda$ be the subdivided claw of $H$ such that
$V(\Lambda) = \{R_0, R_1, R_2, R_3, R_1^+, R_2^+, R_3^+\}$
and $R_0 R_i R_i^+$ is a path of $\Lambda$
for $i = 1,2,3$.
Moreover,
let $R_0 R_i = x_i$ and $R_i R_i^+ = y_i$ for $i=1,2,3$.
Then the following holds.
\begin{itemize}
\item[\normalfont{i)}]
There exists an induced net $N(x_1', x_2', x_3'; y_1',y_2',y_3')$ in $G$
such that
$x_i' \in \{x_i\} \cup l_H(R_0)$ and 
$y_i' \in  \{y_i, x_i\} \cup l_H(R_i) \cup l_H(R_0)$ (see Figure \ref{fig_endmove}).
\item[\normalfont{ii)}]
Each $y_i'$ is heavy. Moreover, if $y_i' \in \{x_i\} \cup l_H(R_0)$, then
$\ed(y_i') \ge \frac{|E(H)|-2}3 + 2 + |J|$,
where $J = \{y_j' \mid j \neq i, \ y_j' \in \{x_j\} \cup l_H(R_0)\}$.
\item[\normalfont{iii)}]
$x_1 x_2 x_3$ is a triangle of $G$ if and only if 
$y_i' \in \{y_i\} \cup l_H(R_i)$ for each $i$.
\end{itemize}
\label{endmove}
\end{cor}
\noindent
\textit{Proof.}\quad
Note that $x_i, y_i \in E(H) = V(G)$.
Since $L(\Lambda)$ is a net,
$\{x_i, y_i \mid 1 \le i \le 3\}$ induces a net
$N(x_1, x_2, x_3; y_1, y_2, y_3)$ in $\cl(G)$.
Moreover, for $0 \le i \le 3$,
$R_i$ corresponds to a clique, say $\hat{R}_i$, in $\cl(G)$ and
$V(\hat{R}_i) \cap \LC(\cl(G))$ corresponds to $l_H(R_i)$.
Thus by Lemma \ref{em_g},
there exists an induced net 
$N' = N(x_1', x_2', x_3'; y_1', y_2', y_3')$ in $G$
such that
$x_i' \in \{x_i\} \cup l_H(R_0)$ and 
$y_i' \in \{y_i, x_i\} \cup l_H(R_i) \cup l_H(R_0)$
for $i=1,2,3$.
Hence i) holds.

Since $N'$ is an induced net of $G$,
$\ed(y_i') = d_{\cl(G)}(y_i') \ge d_G(y_i')
\ge \frac{|V(G)|-2}3 = \frac{|E(H)|-2}3$
for each $i$.
Thus each $y_i'$ is heavy.
Assume that $y_i' \in \{x_i\} \cup l_H(R_0)$ for some $i$.
Without loss of generality, we may assume that $i=1$.
Since $N'$ is an induced net of $G$,
none of the vertex in $\{x_2', x_3'\} \cup J$ is adjacent to $y_1'$ in $G$.
On the other hand,
since $\{x_1, x_2, x_3\} \cup l_H(R_0)$ induces a complete graph in $\cl(G)$
and $ \{x_2', x_3'\} \cup J  \subseteq \{x_2, x_3\} \cup l_H(R_0)$,
each vertex in $\{x_2', x_3'\} \cup J$ is adjacent to $y_1'$ in $\cl(G)$.
Hence
$d_{\cl(G)}(y_1') - d_G(y_1') \ge  2 + |J|$,
and thus
$\ed(y_1') = d_{\cl(G)}(y_1') \ge d_G(y_1') + 2 + |J|
\ge \frac{|V(G)|-2}3 + 2 + |J|
 = \frac{|E(H)|-2}3 + 2 + |J|$.
Therefore ii) holds.
Since $l_H(R_i)$ corresponds to $V(\hat{R}_i) \cap \LC(\cl(G))$,
iii) follows from Lemma \ref{em_g} ii).
\qed
\medskip

A connected multigraph $H$ is called
\textit{essentially $k$-edge-connected}
if $H-F$ has at most one component which contains an edge
for every $F \subseteq E(H)$ with $|F| < k$.
Note that a graph $H$ is essentially $k$-edge-connected
if and only if $L(H)$ is $k$-connected or complete.

\begin{lm}
Let $H$ be an essentially $2$-edge-connected multigraph
and let $x \in V(H)$ such that $d_H(x) \ge 2$.
If $|E(H-\{x\})| \le 3$,
then there exists a DCT of $H$ containing $x$.
\label{nokorisanhen}%
\end{lm}%
\noindent\textit{Proof.}\quad
Let $H$ be the minimal counterexample.
Assume that there exists a cycle $C$ of length $2$ or $3$
containing $x$.
Since $C$ is not a DCT of $H$,
$E(H - V(C)) \neq \emptyset$.
Let $H' = H/C$, then
$H'$ is essentially $2$-edge-connected.
Since $H$ is essentially $2$-edge-connected
and $E(H - V(C)) \neq \emptyset$, we have $d_{H'}(v_C) \ge 2$.
Moreover, the assumption $|E(H-\{x\})| \le 3$ yields
$|E(H'-\{v_C\})| \le 3$.
Hence, by the minimality of $H$,
there exists a DCT of $H'$ containing $v_C$.
Since a cycle of length $2$ or $3$ is collapsible,
by Proposition \ref{col} i),
$H$ has a DCT which contains $V(C)$,
a contradiction.
Hence there exists no cycle of length $2$ or $3$ containing $x$.

Since $H$ does not have a DCT,
$E(H-\{x\}) \neq \emptyset$.
Hence $H$ is not a star with center $x$.
Since $H$ is essentially $2$-edge-connected and $d_H(x) \ge 2$,
there exists a cycle $C$ in $H$ which contains $x$.
If $|V(C)| \ge 5$,
then $|E(C-\{x\})| \ge 3$.
Since $|E(H-\{x\})| \le 3$, we have $E(H-\{x\}) = E(C-\{x\})$. 
Hence $C$ is a DCT of $H$ containing $x$,
a contradiction.
Therefore $|V(C)| = 4$.
Then $|E(C-\{x\})| = 2$ and
there exists a unique edge $z_1z_2$ such that $z_1, z_2 \notin V(C)$.
Since $\{z_1 z_2\}  = E(H-\{x\}) \setminus E(C-\{x\})$,
$N_H(z_1) \cup N_H(z_2) \subseteq \{x, z_1,z_2\}$.
Then we have $xz_1, xz_2 \in E(H)$ or $xz_{i}$ is a multiple edge in $H$ for some $i$, 
since $H$ is essentially $2$-edge-connected.
Hence 
there exists a cycle of length $2$ or $3$
containing $x$,
a contradiction.
\qed
\medskip

Lemma \ref{nokorisanhen} yields the following corollary.
\begin{cor}
Let $H$ be an essentially $2$-connected multigraph
and let $\Xi$ be a collapsible subgraph of $H$.
If $|E(H-\Xi)| \le 3$,
then there exists a DCT of $H$.
\label{nokorisanhen_c}
\end{cor}
\noindent\textit{Proof.}\quad
Let $H' = H / \Xi$,
then $H'$ is an essentially $2$-edge-connected multigraph.
If $v_{\Xi}$ has degree at most $1$ in $H'$,
then $H' \simeq K_1$ or $K_2$
since $H$ is essentially $2$-edge-connected.
Therefore, $v_{\Xi}$ is a DCT of $H'$.
On the other hand,
If $v_{\Xi}$ has degree at least $2$ in $H'$,
then there exists a DCT of $H'$
containing $v_{\Xi}$ by Lemma \ref{nokorisanhen}.
In either case,
$H$ has a DCT by Proposition \ref{col} i).
\qed
\medskip

\begin{figure}[tb]
\begin{center}
\begin{picture}(310,120)
%
%
\includegraphics[width=10cm]{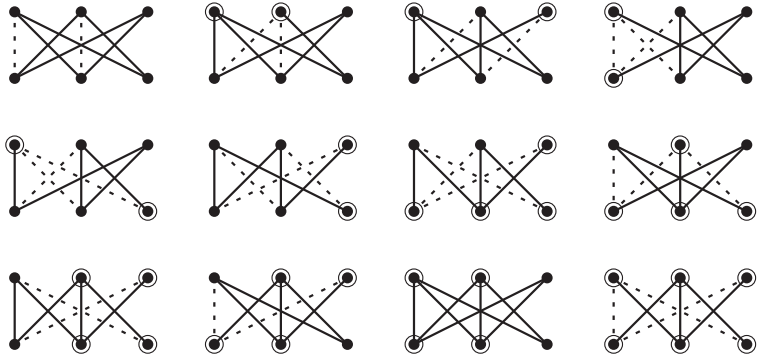}
\end{picture}
\caption{Possible arrangements of $S$ (circled vertices) in $K_{3,3}^-$ and corresponding spanning subgraphs $F$ (solid edges)}
\label{fig_k33-}
\end{center}
\end{figure}

Let $K_{3,3}^-$ be the graph obtained from $K_{3,3}$
by deleting one edge.
By a straightforward case analysis we obtain the following lemma
(see Figure \ref{fig_k33-}).
\begin{lm}
Both of $K_{3,3}$ and $K_{3,3}^-$ are collapsible.
\qed
\label{lm_k33-}%
\end{lm}%
\begin{lm}
Let $H$ be an essentially $2$-edge-connected triangle-free simple graph
which does not contain a DCT
and let $n = |E(H)|$.
If $\{e_1, e_2, e_3\}$ is a matching of $H$,
then $\sum_{i=1}^3 \ed(e_i) \le n+1$.
\label{jisuukeisan}
\end{lm}
\noindent\textit{Proof.}\quad
Assume to the contrary that $\sum_{i=1}^3 \ed(e_i) \ge n+2$.
Let
$\Xi = H[\bigcup_{1 \le i \le 3} V(e_i)]$ and
$\gamma =\sum_{1 \le i < j \le 3}|N_H^e(e_i) \cap N_H^e(e_j)|$.
Then $|E(H-\Xi)|  =  |E(H)| - (|N_H^e(e_1)| + |N_H^e(e_2)| + |N_H^e(e_3)|
 + |\{e_1, e_2, e_3 \}| - \gamma) \le n- (n + 2 + 3 - \gamma) \le \gamma - 5$,
and hence $\gamma \ge 5$.
This yields $|E(\Xi)| \ge 8$.
Since $H$ is triangle-free, $\Xi \simeq K_{3,3}^-$ or $K_{3,3}$,
and hence by Lemma \ref{lm_k33-}, $\Xi$ is collapsible.
Moreover, since $\Xi \simeq K_{3,3}^-$ or $K_{3,3}$,
we have $\gamma \le 6$.
Hence $|E(H-\Xi)| \le \gamma - 5 \le 1$.
By Corollary \ref{nokorisanhen_c},
$H$ has a DCT,
a contradiction.
\qed
%
%
%
%
%
%
%
%
%
%
%
%
%
%
%
%
%
\section{Proof for the large graphs}
\label{large}
\begin{figure}[tb]
\begin{center}
\begin{picture}(200,120)
\put(-10,65){\footnotesize$A$}
\put(10,105){\footnotesize$a_1$}
\put(30,76){\footnotesize$a_2$}
\put(35,45){\footnotesize$a_3$}
\put(20,135){\footnotesize$C_{1,1}$}
\put(50,135){\footnotesize$C_{1,2}$}
\put(78,135){\footnotesize$C_{1,3}$}
\put(100,135){\footnotesize$C_{1,p_1-1}$}
\put(135,135){\footnotesize$C_{1,p_1}$}
\put(37,120){\footnotesize$c_{1,1}$}
\put(64,120){\footnotesize$c_{1,2}$}
\put(124,118){\footnotesize$b_1'$}
\put(50,60){\footnotesize$C_{2,1}$}
\put(110,60){\footnotesize$C_{2,2}$}
\put(70,76){\footnotesize$c_{2,1}=b_2'$}
\put(78,40){\footnotesize$C_{3,1}$}
\put(88,15){\footnotesize$c_{3,1}=b_3'$}
\put(88,-5){\footnotesize$D_3$}
\put(172,65){\footnotesize$B$}
\put(152,105){\footnotesize$b_1$}
\put(140,76){\footnotesize$b_2$}
\put(133,45){\footnotesize$b_3$}
\includegraphics[width=6cm]{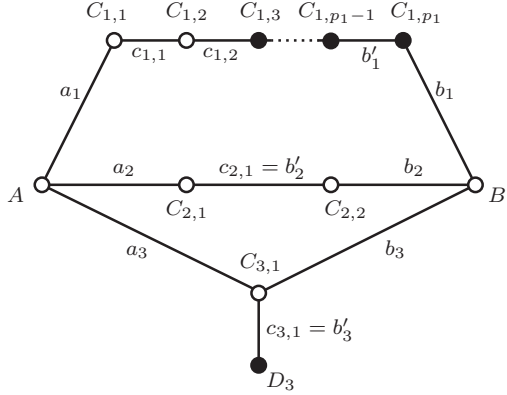}
\end{picture}
\caption{The graph $\T$ with $p_2 = 2$ and $p_3 = 1$,
where the white vertices denote $X$.}
\label{theta}
\end{center}
\end{figure}

In this section we prove Conjecture \ref{conj}
for the graphs with at least $33$ vertices.
As we will see at the end of this section,
the proof immediately follows from the following theorem.
We call a matching \textit{heavy}
if each edge of the matching is heavy.
\begin{thm}
Let $G$ be a graph of order $n$
which satisfies the assumption of Conjecture \ref{conj}
and let $H$ be the triangle-free graph such that
$L(H) = \cl(G)$.
Then there exists either a DCT
or a heavy matching of size $4$ in $H$.
\label{main}
\end{thm}
\noindent
\textit{Proof.}\quad
Suppose that $H$ has neither a DCT
nor a heavy matching of size $4$.
Then,
by Theorem \ref{HN},
$\cl(G)$ is not hamiltonian,
and hence it follows from Theorem \ref{Br} that
$\cl(G)$ contains an induced subgraph $F \in \F$.
Let $\T$ be the subgraph of $H$ such that
$L(\T) = F$.
Then there exist two vertices $A$, $B$ of degree $3$ and
three internally vertex-disjoint paths $P_i$ ($1 \le i \le 3$)
of length at least two joining $A$ and $B$ in $\T$.
Moreover, if $|P_i| = 3$, 
then the middle vertex of $P_i$ is joined to one pendant edge in $\T$
(that is, the middle vertex of $P_i$ has degree $3$ in $\T$).
We denote
$p_i = |P_i| - 2$.
Let
$V(\T) \setminus \{A,B\} = 
\{C_{i,j} \mid 1 \le i \le 3, 1 \le j \le p_i\} \cup
\{D_i \mid 1 \le i \le 3, p_i=1\}$,
where $P_i = A C_{i,1} C_{i,2} \ldots C_{i,p_i} B$
for $1 \le i \le 3$
and $C_{i,1} D_i$ is the pendant edge
if $p_i = 1$.
For $1 \le i \le 3$,
let $a_i = AC_{i,1}$ and $b_i = C_{i,p_i} B$.
Moreover,
if $p_i \neq 1$, then
let $c_{i,j} = C_{i,j} C_{i,j+1}$ for $1 \le j \le p_i -1$,
and  if $p_i = 1$, then
let $c_{i,1} = C_{i,1} D_i$.
Since $F=L(\T)$ is an induced subgraph of $\cl(G)$,
$D_i \neq D_j$ and $D_i \notin V(P_j)$ hold for $i \neq j$.
Let $X =(\{C_{i,j} \mid 1 \le i \le 3, j=1,2\} \cap V(H)) \cup \{A,B\}$
(see Figure \ref{theta}).

For $1 \le i \le 3$,
let $b_i' = C_{i,p_i-1} C_{i,p_i}$ if $p_i \neq 1$
and let
$b_i' = c_{i,1}$ if $p_i = 1$.
Since $\{a_i, c_{i,1} \mid 1 \le i \le 3\}$ and $\{b_i, b_i' \mid 1 \le i \le 3\}$
induces a subdivided claw in $H$,
it follows from Corollary \ref{endmove} i) that
there exists a heavy edge $c_i^*$ in $\{c_{i,1}, a_i\} \cup l(C_{i,1}) \cup l(A)$ and
a heavy edge $\bar{c}_i^*$ in $\{b_i, b_i'\} \cup l(C_{i,p_i}) \cup l(B)$
for each $i$.
\medskip
\\
%
%
%
%
%
%
%
\noindent
\textbf{Case 1.}\quad
$c_i^* \in \{c_{i,1}\} \cup l(C_{i,1})$ holds for $i=1,2,3$ or
$\bar{c}_i^* \in \{b_i'\} \cup l(C_{i,p_i})$ holds for $i=1,2,3$.
\medskip

Without loss of generality, we may assume that
$c_i^* \in \{c_{i,1}\} \cup l(C_{i,1})$ holds for $i=1,2,3$.
Then by Corollary \ref{endmove} iii),
$a_1 a_2 a_3$ is a triangle in $G$.
\begin{claim}\ 
\begin{itemize}
\item[\normalfont{i)}]
For each $i$ with $p_i = 1$,
there exists
$b_i^* \in \{b_i\} \cup l(C_{i,1})$
such that $b_i^*$ is heavy in $H$.%
\item[\normalfont{ii)}]
If there exists
$u \in E(H - X)$
such that $u$ and $c_i^*$ are adjacent in $H$,
then $p_i = 1$ and $c_i^* = c_{i,1}$.
\end{itemize}
\label{ichinisan}
\end{claim}
\noindent
\textit{Proof.}\quad
If $p_i = 1$, then
$\{a_i, b_i\} \cup \{a_j, c_{j,1} \mid 1 \le j \le 3, j\neq i\}$
induces a subdivided claw in $H$.
Hence by Corollary \ref{endmove} i) and ii)
there exists a heavy edge $b_i^*$ in $\{b_i, a_i\} \cup l(C_{i,1}) \cup l(A)$.
Since $a_1 a_2 a_3$ is a triangle in $G$,
by Corollary \ref{endmove} iii),
we have $b_i^* \in \{b_i\} \cup l(C_{i,1})$.
Therefore i) holds.

If there exists $u$ as in ii),
then $u$ is incident with $D_i$,
since $u$ is incident with neither
$C_{i,1}$ nor $C_{i,2}$.
Hence $p_i = 1$ and $c_i^* = C_{i,1}D_i = c_{i,1}$.
\qed
\begin{claim}
$p_i \le 3$
for any $i$ with $1 \le i \le 3$.
\label{hachi}
\end{claim}
\noindent\textit{Proof.}\quad
Assume to the contrary that $p_i \ge 4$ for some $i$.
Recall that
there exists a heavy edge $\bar{c}_i^* \in \{b_i', b_i\} \cup l(C_{1,p_1}) \cup l(B)$.
For each $j$ with $1 \le j \le 3$,
$\bar{c}_i^*$ and $c_j^*$ are not adjacent in $H$
since $c_j^* \in \{c_{j,1} \} \cup l(C_{i,1})$.
Therefore $\{\bar{c}_i^*, c_1^*, c_2^*, c_3^*\}$ is a heavy matching of size $4$,
a contradiction.
\qed
\begin{claim}
Let $UW \in E(H)$ such that
$U$ is adjacent to $A$ or $B$.
Then $UW$ is dominated by $X$.
\label{nana}
\end{claim}
\noindent\textit{Proof.}\quad
Assume to the contrary that
the edge $u=UW$ is not dominated by $X$.
Then
$U,W \notin \{C_{i,1}, C_{i,2}, A, B\}$ for any $i$.
Let $R = A$ if $U$ is adjacent to $A$,
and let $R=B$ otherwise.
Moreover, let $r$ be the edge of $H$ joining $U$ and $R$.
We shall prove that there exists a subdivided claw
containing both $u$ and $r$ with center $R$.

Since $D_l$ ($p_l=1$) are distinct vertices, 
we can take $i,j$ so that $W \neq D_i, D_j$.
Moreover, since $H$ is triangle-free,
we have $U \neq D_l$ for each $l$.
Hence, if $U \in N_H(A)$, then
$\{r, u, a_i, c_{i,1}, a_j, c_{j,1}\}$ is a desired subdivided claw.
Thus we assume $U \in N_H(B)$.
Let $k = \{1,2,3\} \setminus \{i,j\}$.
If $U = C_{i,3}$,
either $\{r, u, b_j, b_j', b_k, b_k'\}$ (in the case $W \neq D_k$) or
$\{r, u, b_j, b_j', b_k, a_k\}$ (in the case $W = D_k$)
induces a desired subdivided claw.
The case $U = C_{j,3}$ is similar,
and hence we may suppose $U \neq C_{i,3}, C_{j,3}$.
Since $H$ is triangle-free,
we have $W \neq C_{l,3}$ for each $l$.
Hence $u$ and $b_l$ are not adjacent in $H$ for $l=i,j$ and,
since $u$ is not dominated by $X$,
$u$ and $b_l'$ are not adjacent in $H$ for $l=i,j$.
Therefore
$\{b_i, b_i', b_j, b_j', r, u\}$ induces a desired subdivided claw.

By Corollary \ref{endmove} i) and ii),
there exists $u^* \in \{u,r\} \cup l(U) \cup l(R)$ which is heavy.
Since there exists no heavy matching of size $4$ in $H$,
we may assume that $u^*$ and $c_1^*$ are adjacent.
Since both of the endvertices of $c_1^*$ are contained in
$\{C_{1,1}, C_{1,2}, D_1\} \cup V_1(H)$
and
$U,W \notin \{C_{1,1}, C_{1,2}\}$,
we have $W=D_1$,
which implies $p_1 = 1$.
Then by Claim \ref{ichinisan} i),
there exists $b_1^* \in \{b_1\} \cup l(C_{1,1})$ which is heavy.
Since both of the endvertices of $b_1^*$ are contained in
$\{C_{1,1}, B\} \cup V_1(H)$ and
both of the endvertices of $c_i^*$ are contained in
$\{C_{i,1}, C_{i,2}, D_i\} \cup V_1(H)$ for $i=2,3$,
$\{b_1^*, c_2^*, c_3^*, u^*\}$
is a heavy matching of size $4$, a contradiction.
\qed
%
%
%
%
%
%
%
%
\begin{claim}
For any $i$ with $p_i = 3$,
$N_H(C_{i,3}) \subseteq X$.
\label{March7_maru2}
\end{claim}
\noindent\textit{Proof.}\quad
Assume to the contrary that
there exists $U \in N_H(C_{i,3}) \setminus X$.
Without loss of generality,
we may assume that $i=1$.
Let $u = U C_{1,3}$.
Since $b_2'$ and $b_3'$ have no common endvertex,
we may assume without loss of generality that
$b_2'$ is not incident with $U$.
Note that $U \neq C_{3,3}$ in the case $p_3 = 3$,
since $H$ is triangle-free.
Hence, if $b_3'$ is incident with $U$,
then $p_3 = 1$ and $U=D_3$,
because $U \notin X$.
Therefore,
if $p_3 \neq 1$, then
$\{b_1, u, b_2, b_2', b_3, b_3'\}$ induces a subdivided claw,
and
if $p_3 = 1$, then
$\{b_1, u, b_2, b_2', b_3, a_3\}$ induces a subdivided claw.
In either case,
by Corollary \ref{endmove} i) and ii),
there exists a heavy edge $b_1^{**} \in \{u, b_1\} \cup l(C_{1,3}) \cup l(B)$.
Since $U \notin X$,
$u$ and $c_i^*$ are not adjacent for $i=1,2$.
Hence we can deduce that
$\{b_1^{**}, c_1^*, c_2^*\}$ is a heavy matching of $H$.
Since there exists no heavy matching of size $4$,
$b_1^{**}$ and $c_3^*$ are adjacent.
This implies that
$p_3 = 1$, $b_1^{**} = u$ and $U = D_3$.
Hence
$\{b_1, u, b_2, b_2', b_3, a_3\}$ induces a subdivided claw.
By Corollary \ref{endmove} i) and ii),
there exists a heavy edge $a_3^{**} \in \{a_3, b_3\} \cup l(C_{3,1}) \cup l(B)$.
Now $b_1^{**} = u = D_3 C_{1,3}$ yields
$\{b_1^{**}, c_1^*, c_2^*, a_3^{**}\}$ is a heavy matching of $H$,
a contradiction.
\qed
\begin{claim}
If $u \in E(H)$ is not dominated by $X$,
then $u$ is not heavy in $H$.
\label{santengo}
\end{claim}
\noindent
\textit{Proof.}\quad
Assume, to the contrary, that $u$ is heavy.
Since $H$ does not contain a heavy matching of size $4$,
$u$ and $c_i^*$ are adjacent in $H$ for some $i$.
Without loss of generality, we may assume that $i=1$.
Then by Claim \ref{ichinisan} ii) we have
$p_1 = 1$ and $c_1^* = c_{1,1}$,
and then
we can take $b_1^*$ as in Claim \ref{ichinisan} i).
Since both of the endvertices of $b_1^*$ are contained in $\{C_{1,1}, B\} \cup V_1(H)$,
$u$ and $b_1^*$ are not adjacent in $H$.

Since $\{b_1^*, c_2^*, c_3^*\}$ is a heavy matching of $H$
and there exists no heavy matching of size $4$ in $H$,
$u$ and $c_i^*$ are adjacent in $H$ for $i=2$ or $3$.
Without loss of generality, we may assume that
$u$ and $c_2^*$ are adjacent.
Then by Claim \ref{ichinisan} ii),
$p_2 = 1$ and $c_2^* = c_{2,1}$,
and hence there exists $b_2^*$ as in Claim \ref{ichinisan} i).

Recall that both of $u, c_{1,1}$ and $u, c_{2,1}$ are adjacent in $H$.
Since $u$ is not dominated by $X$,
we have $u = D_1 D_2$.
Thus $\{a_1, a_2, c_1^*, u, b_1, b_3\}$ induces a subdivided claw in $H$.
By Corollary \ref{endmove} i) and ii),
there exists $a^* \in \{a_2, a_1\} \cup l(A) \cup l(C_{1,1})$ which is heavy.
If $a^* \in \{a_1\} \cup l(C_{1,1})$,
then $\{a^*, u, b_2^*, c_3^*\}$ is a heavy matching of $H$,
a contradiction.
On the other hand,
if $a^* \in \{a_2\} \cup l(A)$,
then $\{a^*, u, b_1^*, c_3^*\}$ is a heavy matching of $H$,
a contradiction.
\qed
\medskip

Let $\Q_1$ be the set of paths in $H-E(P_1 \cup P_2 \cup P_3)$
joining $C_{i,1}$ and $C_{j,1}$ for some $i,j$ with $i \neq j$
and $p_i = p_j = 1$,
and let 
$\Q_2$ be the set of paths in $H-E(P_1 \cup P_2 \cup P_3)$
joining $C_{i,1}$ and $C_{j,2}$ ($i \neq j$) or
$C_{i,1}$ and $C_{j,1}$ ($i \neq j$, $p_i = 1$ and $p_j \ge 2$).
Note that a path in $\Q_1 \cup \Q_2$ may contain the edge $C_{i,1} D_i$
for some $i$ with $p_i=1$.
\begin{claim}
$\Q_1 \cup \Q_2 \neq \emptyset$.
\label{yon}
\end{claim}
\noindent
\textit{Proof.}\quad
Assume, to the contrary, that $\Q_1 = \Q_2 = \emptyset$.
Let $i,j \in \{1,2,3\}$ with $i \neq j$.
If $p_i=p_j=1$,
then $C_{i,1}D_j, D_i D_j \notin E(H)$ since $\Q_1 = \emptyset$.
If $p_i=1$ and $p_j \ge 2$,
then $C_{i,1} C_{j,2}, D_i C_{j,1}, D_i C_{j,2} \notin E(H)$ since $\Q_2 = \emptyset$.
If $p_i, p_j \ge 2$,
then $C_{i,1} C_{j,2} \notin E(H)$ since $\Q_2 = \emptyset$.
Moreover, in either case, $C_{i,1} C_{j,1} \notin E(H)$ since $H$ is triangle-free.

By the above argument,
it follows that
$N_H^e(c_i^*) \cap N_H^e(c_j^*) \subseteq \{C_{i,2} C_{j,2}\}$
for any $i,j$ with $i \neq j$.
Note that if $p_i \ge 3$, then $b_i \notin N_H^e(c_j^*)$ for any $j$.
Let $E_0 = E(H) \setminus
( N_H^e(c_1^*) \cup N_H^e(c_2^*)  \cup N_H^e(c_3^*)
\cup \{c_1^*, c_2^*, c_3^*\}$
$\cup\; \{b_i \mid p_i \ge 3\})$,
then
\begin{eqnarray*}
n & = & |E(H)| \\
& \ge &
\left| N_H^e(c_1^*) \cup N_H^e(c_2^*)  \cup N_H^e(c_3^*) |
+ | \{c_1^*, c_2^*, c_3^*\}  \right|
 + |\{b_i \mid p_i \ge 3\}| + |E_0|
\\
& \ge & |N_H^e(c_1^*)| + |N_H^e(c_2^*)| + |N_H^e(c_3^*)|
 - |\{C_{i,2} C_{j,2} \mid 1 \le i < j \le 3\} \cap E(H)| + 3 
 + |\{b_i \mid p_i \ge 3\}| + |E_0| \\
& \ge & 3 \times \frac{n-2}3 + 3 - |\{C_{i,2} C_{j,2} \mid 1 \le i < j \le 3\} \cap E(H)|
 + |\{b_i \mid p_i \ge 3\}| + |E_0|,
\end{eqnarray*}
and hence
\begin{eqnarray}
|\{b_i \mid p_i \ge 3\}| + |E_0| 
\le |\{C_{i,2} C_{j,2} \mid 1 \le i < j \le 3\} \cap E(H)| -1.
\label{eq_qexist}
\end{eqnarray}
Without loss of generality,
we may assume that $p_1 \ge p_2 \ge p_3$.
Let $t= |\{i \mid p_i \ge 3\}| = |\{b_i \mid p_i \ge 3\}|$.
If $t=0$, then
$|\{C_{i,2} C_{j,2} \mid 1 \le i < j \le 3\} \cap E(H)|=0$
since $H$ is triangle-free.
Then the right hand side of (\ref{eq_qexist}) is $-1$,
a contradiction.
If $t \ge 2$, then by (\ref{eq_qexist}),
$|\{C_{i,2} C_{j,2} \mid 1 \le i < j \le 3\} \cap E(H)| \ge 3$.
This implies that $C_{1,2} C_{2,2} C_{3,2}$ is a triangle of $H$,
a contradiction.
Hence $t = 1$, which yields $p_1 = 3$ and $p_2, p_3 \le 2$.
Then we have $C_{2,2}C_{3,2} \notin E(H)$,
since otherwise $C_{2,2}C_{3,2}B$ is a triangle of $H$.
Hence, by (\ref{eq_qexist}),
$E_0 = \emptyset$ and
$C_{1,2}C_{2,2}, C_{1,2}C_{3,2} \in E(H)$.
This yields $p_2 = p_3 = 2$ and $c_i^* \in \{C_{i,1} C_{i,2}\} \cup l(C_{i,1})$ for each $i$.
In the case $C_{1,3} C_{2,1} \in E(H)$,
let $T =A C_{1,1} C_{1,2} C_{1,3} C_{2,1} C_{2,2} B C_{3,2} C_{3,1} A$.
Then,
since $T$ contains $\{C_{i,1}, C_{i,2}\}$ for each $i$,
$T$ dominates every edge of $N_H^e(c_i^*) \cup \{c_i^*\}$.
Moreover, since $V(b_1) \subseteq V(T)$ and $E_0 = \emptyset$, $T$ dominates $E(H)$, a contradiction.
Hence we have $C_{1,3} C_{2,1} \notin E(H)$.
By symmetry we have $C_{1,3} C_{3,1} \notin E(H)$.

Recall that there exists a heavy edge $\bar{c}_1^* \in \{b_1', b_1\} \cup l(C_{1,3}) \cup l(B)$.
If $\bar{c}_1^* \neq b_1'$, then
$\{c_1^*, c_2^*, c_3^*, \bar{c}_1^*\}$ is a heavy matching of $H$,
a contradiction.
Hence $\bar{c}_1^* = b_1'$.
Since $H$ is triangle-free and $p_2=p_3=2$,
we have
$C_{1,3} C_{i,2} \notin E(H)$ for $i=2,3$.
Hence $N_H^e(\bar{c}_1^*) \cap N_H^e(c_i^*) \subseteq \{C_{1,2} C_{i,2}\}$ for $i=2,3$.
Moreover, since $C_{2,2} C_{3,2} \notin E(H)$,
$N_H^e(c_2^*) \cap N_H^e(c_3^*) = \emptyset$.

Let $E_1 = E(H) \setminus
\left(N_H^e(\bar{c}_1^*) \cup N_H^e(c_2^*)  \cup N_H^e(c_3^*)
\cup \{\bar{c}_1^*, c_2^*, c_3^*\}\right)$,
then
\begin{eqnarray*}
n & = & |E(H)| \\
& \ge &
\left| N_H^e(\bar{c}_1^*) \cup N_H^e(c_2^*)  \cup N_H^e(c_3^*) |
+ | \{\bar{c}_1^*, c_2^*, c_3^*\}  \right| + |E_1|
\\
& \ge & |N_H^e(\bar{c}_1^*)| + |N_H^e(c_2^*)| + |N_H^e(c_3^*)|
 - |\{C_{1,2} C_{2,2} , C_{1,2} C_{3,2} \}| + 3 
 + |\{a_1\}|\\
& \ge & 3 \times \frac{n-2}3 - 2 + 3 + 1\, =\, n,
\end{eqnarray*}
and hence $E_1 = \{a_1\}$.
Let $T' = A C_{2,1} C_{2,2} C_{1,2} C_{1,3} B C_{3,2} C_{3,1} A$.
Then, since
$\{C_{1,2}, C_{1,3}, C_{2,1}, C_{2,2}, C_{3,1}, C_{3,2}\} \subseteq V(T')$,
$T'$ dominates every edge of $N_H^e(\bar{c}_1^*) \cup N_H^e(c_2^*)  \cup N_H^e(c_3^*)
\cup \{\bar{c}_1^*, c_2^*, c_3^*\}$.
Since $a_1$ is dominated by $T'$,
$T'$ is a DCT of $H$, a contradiction.
\qed
\medskip

In the case where $\Q_1 \neq \emptyset$,
take $Q_1, Q_2, \ldots , Q_m \in \Q_1$ so that
$|V(Q_1)| + \ldots + |V(Q_m)|$ is as large as possible,
subject to the condition that $Q_1, \ldots , Q_m$ are
internally vertex-disjoint,
and let $\Q = \{Q_1, Q_2, \ldots , Q_m\}$.
In the case where $\Q_1 = \emptyset$,
take $Q \in \Q_2$
and let $\Q = \{Q\}$.
\begin{claim}
There exists a closed trail $T$ of $H$ such that
$X \cup V(P_i) \cup V(P_j) \cup
\left( \bigcup_{Q \in \Q} V(Q)\right) \subseteq V(T)
\subseteq \left( \bigcup_{l=1}^3 V(P_l) \right) \cup  \left( \bigcup_{Q \in \Q} V(Q) \right)$
for some $i,j$ with $1 \le i<j \le 3$.
\label{yontengo}
\end{claim}
\noindent
\textit{Proof.}\quad
Assume $\Q_1 = \emptyset$ and let $Q$ be the (unique) path in $\Q$.
Then without loss of generality,
we may assume that
either
$Q$ joins $C_{1,1}$ and $C_{2,2}$
or 
$p_1 = 1$, $p_2 \ge 2$ and $Q$ joins $C_{1,1}$ and $C_{2,1}$.
In the former (resp.~latter) case,
$A C_{2,1} C_{2,2} Q C_{1,1} C_{1,2} \ldots C_{1,p_1} B P_3 A$
(resp.~$A C_{1,1} Q C_{2,1} C_{2,2} \ldots C_{2,p_2} B P_3 A$)
is a required closed trail, where $i = 1$ and $j=3$.
Hence we may assume that $\Q_1 \neq \emptyset$.

We apply induction on $|\Q|$,
and we find the desired closed trail without using the assumption that
$H$ is essentially $2$-edge-connected.
In the case $|\Q| = 1$,
we may assume without loss of generality that $Q_1$ joins $C_{1,1}$ and $C_{2,1}$
and $p_1 = p_2 = 1$. Then $A C_{1,1} Q C_{2,1} B P_3 A$ is a required closed trail.
Suppose that $|\Q| = 2$.
If
$Q_1$ and $Q_2$ have the same endvertices,
say $C_{1,1}$ and $C_{2,1}$,
then
$A C_{1,1} Q_1 C_{2,1} Q_2 C_{1,1} B P_3 A$ is a required closed trail.
Otherwise, without loss of generality we may assume that
$Q_1$ joins $C_{1,1}$ and $C_{2,1}$ and 
$Q_2$ joins $C_{2,1}$ and $C_{3,1}$.
Then
$A C_{1,1} Q_1 C_{2,1} B C_{3,1} Q_2 C_{2,1} A$ is a required closed trail.

Assume that $|\Q| \ge 3$.
If $|\Q| = 3$ and
$Q_1$ joins $C_{i,1}$ and $C_{j,1}$,
$Q_2$ joins $C_{j,1}$ and $C_{k,1}$ and
$Q_3$ joins $C_{k,1}$ and $C_{i,1}$ for some $i,j,k$ with $\{i,j,k\} = \{1,2,3\}$,
then
$A C_{i,1} B C_{j,1} Q_2 C_{k,1} Q_3 C_{i,1} Q_1 C_{j,1} A$
is a required closed trail.
Otherwise, there exist $Q_a, Q_b \in \Q$ such that
$Q_a$ and $Q_b$ have the same endvertices.
Let $\Q' = \Q \setminus \{Q_a, Q_b\}$,
then by the induction hypothesis,
there exists a closed trail $T$ in $H-E(Q_a \cup Q_b)$ such that
$X \cup V(P_i) \cup V(P_j) \cup
\left( \bigcup_{Q \in \Q'} V(Q)\right) \subseteq V(T)
\subseteq \left( \bigcup_{l=1}^3 V(P_l) \right) \cup  \left( \bigcup_{Q \in \Q'} V(Q) \right)$
for some $i,j$,
and then $T \cup Q_a \cup Q_b$ is a required closed trail.
\qed
\begin{claim}
Let $P_u = RUW$ be the path of length two in $H$
such that $R \in V(T)$ and the edge $UW$ is not dominated by $T$.
Then
$R =C_{i,1}$ or $C_{i,2}$ for some $i$.
\label{March7_maru1_maru4}
\end{claim}
\noindent\textit{Proof.}\quad
Assume to the contrary that
$R \neq C_{i,1},C_{i,2}$ for any $i$.
By Claim \ref{nana}, we have $R \neq A,B$,
and by Claim \ref{March7_maru2}, we have
$R \neq C_{i,3}$ for any $i$.
Hence $R$ is an internal vertex
of a path $Q \in \Q$.
Without loss of generality,
we may assume that $Q$ joins
either $C_{1,1}$ and $C_{2,1}$ or $C_{1,1}$ and $C_{2,2}$.

Let $Q^1$ be the path in $P_1 \cup Q$ which joins $R$ and $B$,
and
let $Q^2$ be the path in $P_2 \cup Q$ which joins $R$ and $A$.
Moreover, let $\tilde{Q}^1$ (resp.~$\tilde{Q}^2$) be the subpath of $Q^1$ (resp.~$Q^2$)
of length two which contains $R$.
Then both of $\tilde{Q}^1$ and $\tilde{Q}^2$
are contained in $P_1 \cup P_2 \cup Q$.
Since $UW$ is not dominated by $T$,
$E(\tilde{Q}^1 \cup \tilde{Q}^2 \cup P_u)$ induces a subdivided claw.
By Corollary \ref{endmove} i) and ii),
there exists a heavy edge $u^* \in \{WU, UR\} \cup l(U) \cup l(R)$.
Since $R \notin X$, $u^*$ is not dominated by $X$,
which contradicts Claim \ref{santengo}.
\qed
\begin{claim}
Let $v_a, v_b \in E(H)$ such that
$v_a \in \{a_i\} \cup l(C_{i,1}) \cup l(A)$ and
$v_b \in \{b_i\} \cup l(C_{i,p_i}) \cup l(B)$
for some $i$ with $p_i \ge 2$.
Then $\{v_a, v_b\}$ is not a heavy matching.
\label{March7_maru4ue}
\end{claim}
\noindent\textit{Proof.}\quad
If $\{v_a, v_b\}$ is a heavy matching,
then $\{v_a, v_b, c_j^*, c_k^*\}$ is a heavy matching,
where $\{j,k\} = \{1,2,3\} \setminus \{i\}$.
This is a contradiction.
\qed
\medskip

\begin{figure}[tb]
\begin{center}
\begin{picture}(200,40)
\put(-10,37){\footnotesize$T$}
\put(175,37){\footnotesize$T$}
\put(3,20){\footnotesize$R^1$}
\put(28,26){\footnotesize$r_1$}
\put(59,27){\footnotesize$s_1$}
\put(42,9){\footnotesize$S^1$}
\put(70,-10){\footnotesize$Q_u^1$}
\put(133,26){\footnotesize$r_2$}
\put(118,9){\footnotesize$S^2$}
\put(131,-10){\footnotesize$Q_u^2$}
\put(93,39){\footnotesize$u=s_2$}
\put(157,20){\footnotesize$R^2$}
\includegraphics[width=6cm]{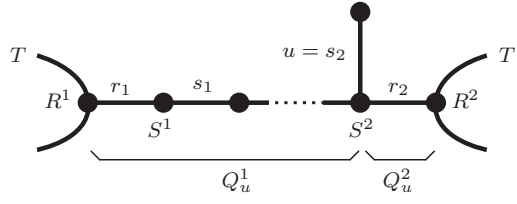}
\end{picture}
\caption{The case where $Q_u^1 \cap \Q_u^2 \neq \emptyset$ and $R^1 \neq R^2$.}
\label{fig_QRS}
\end{center}
\end{figure}

Since $H$ does not have a DCT,
there exists $u \in E(H)$ which is not dominated by $T$.
Since $H$ is essentially $2$-edge-connected, 
there exist two edge-disjoint paths
$Q_u^1$, $Q_u^2$ each of which joins
an endvertex of $u$ and a vertex in $T$.
For $i=1,2$,
take such $Q_u^i$ so that
$|V(Q_u^i) \cap V(T)| = 1$,
and let $R^i$ be the vertex in $V(T) \cap V(Q_u^i)$
(see Figure \ref{fig_QRS}).
Let $r_i$ be the edge of $Q_u^i$
which is incident with $R^i$,
then we can take the edge $s_i \in E(Q_u^i) \cup \{u\}$
so that $s_i$ and  $r_i$ are adjacent in $H$
and $s_i$ is not dominated by $T$.
Hence it follows from Claim \ref{March7_maru1_maru4} that
\begin{eqnarray}
\mbox{$R^i = C_{j,1}$ or $C_{j,2}$ for some $j \in \{1,2,3\}$.}
\label{Qnotanntenn}
\end{eqnarray}
By Claim \ref{March7_maru2},
$C_{j,3} \notin V(Q_u^i)$ for any $j$ with $p_j = 3$.
Since $s_i$ is not dominated by $X$,
$s_i$ is not dominated by $\bigcup_{l=1}^3 P_l$ as well.

In the case $R^1 \neq R^2$,
let $Q_u$ be the path joining $R^1$ and $R^2$
which is contained in $Q_u^1 \cup Q_u^2 \cup \{u\}$,
and in the case  $R^1 = R^2$,
let $Q_u$ be the maximal closed trail
which is contained in $Q_u^1 \cup Q_u^2 \cup \{u\}$.
Moreover, let $S^i$ be the common endvertex of the two edges $r_i$ and $s_i$.
\begin{claim}
For $i=1,2$,
there exists $r_i^*$ such that
$r_i^* \in \{r_i\} \cup l(R^i)$ and
$\ed(r_i^*) \ge \frac{n-2}3 + 2$.
\label{siheavy}
\end{claim}
\noindent\textit{Proof.}\quad
Let $i \in \{1,2\}$.
By (\ref{Qnotanntenn}),
we may assume $R^i \in V(P_1)$ without loss of generality.
Recall that $s_i$ is not dominated by $\bigcup_{l=1}^3 P_l$.
If $R^i=C_{1,1}$ and $p_1 = 1$,
then there exists a subdivided claw
induced by $\{r_i, s_i, a_1, a_2, b_1, b_3\}$.
If $R^i=C_{1,1}$ and $p_1 \ge 2$,
then there exists a subdivided claw
induced by $\{r_i, s_i, a_1, a_2\}$ and two edges in $E(P_1) \setminus \{a_1\}$.
If $R^i=C_{1,2}$,
then there exists a subdivided claw
induced by $\{r_i, s_i, c_{1,1}, a_1, b_1, b_2\}$
or $\{r_i, s_i, c_{1,1}, a_1, c_{1,2}, b_1\}$.
In either case,
the assertion follows from Corollary \ref{endmove} i) and ii),
since Claim \ref{santengo} implies that
each edge in $\{s_i\} \cup l(S^i)$ is not heavy.
\qed
\begin{claim}
$R^1 = R^2$.
\label{hachigatu_new}
\end{claim}
\noindent\textit{Proof.}\quad
By (\ref{Qnotanntenn}), $R^1 \in \{C_{i,1}, C_{i,2}\}$
and $R^2 \in \{C_{j,1}, C_{j,2}\}$ for some $i, j \in \{1,2,3\}$.
Recall that $V(Q_u) \cap V(T) = \{R^1, R^2\}$.
Hence $Q_u$ and any path in $\Q$ are internally vertex-disjoint.
Suppose that $p_i = p_j = 1$.
If $i \neq j$, then $Q_u \in \Q_1$,
which contradicts the maximality of $\Q$.
On the other hand, if $i=j$, then
we have $R^1 = R^2 = C_{i,1}$,
and hence the assertion holds.
Thus it suffices to consider the case $\max\{p_i, p_j\} \ge 2$.

Without loss of generality, we may assume that
$R^1 \in \{C_{1,1}, C_{1,2}\}$, $p_1 \ge 2$
and $R^2 \in \{C_{1,1}, C_{1,2}, C_{2,1}, C_{2,2}\}$
with $R^1 \neq R^2$.
Take $r_1^*$ and $r_2^*$ as in Claim \ref{siheavy}.
Suppose that $r_1^*$ and $r_2^*$ are not adjacent,
then, if $p_3 \ge 2$,
 $\{r_1^*, r_2^*, c_3^*\}$ is a matching with
$\ed(r_1^*) + \ed(r_2^*) + \ed(c_3^*) \ge n+2$, and
if $p_3 = 1$, $\{r_1^*, r_2^*, b_3^*\}$ is a matching with
$\ed(r_1^*) + \ed(r_2^*) + \ed(b_3^*) \ge n+2$.
This contradicts Lemma \ref{jisuukeisan},
and hence $r_1^*$ and $r_2^*$ are adjacent.
Since $R^1 \neq R^2$,
we have $S^1 = S^2$.
If $R^2 \in \{C_{1,1}, C_{1,2}\}$,
then $R^1 S^1 R^2$ is a triangle of $H$,
a contradiction.
Hence we may assume that $R^2 \in \{C_{2,1}, C_{2,2}\}$.

Suppose $S^1 = S^2 = D_j$ for some $j\in \{2,3\}$ with $p_j = 1$.
Then it follows from the fact $p_1 \ge 2$ that
$\{c_{j,1}, s_1, a_j, a_1, b_j, b_1\}$ induces a subdivided claw.
By Corollary \ref{endmove} i) and ii),
there exist heavy edges $a^*$ and $b^*$ such that
$a^* \in \{a_1, a_j\} \cup l(A) \cup l(C_{j,1})$
and $b^* \in \{b_1, b_j\} \cup l(B) \cup l(C_{j,1})$.
By Claim \ref{March7_maru4ue}
we have either $a^* \notin \{a_1\} \cup l(A)$
 or $b^* \notin \{b_1\} \cup l(B)$,
and hence $a^* \in \{a_j\} \cup l(C_{j,1})$
or $b^* \in \{b_j\} \cup l(C_{j,1})$.
Let $c_j^{**} = a^*$ (resp.~$b^*$)
in the former (resp.~latter) case,
then
it follows from Corollary \ref{endmove} i) and ii) that
$\ed(c_j^{**}) \ge \frac{n-2}3 + 2$.
Hence $\{r_1^*, c_j^{**}, c_k^* \}$ is a matching with
$\ed(r_1^*) + \ed(c_j^{**}) + \ed(c_k^*) \ge n+2$,
where $k \in \{2,3\}\setminus \{j\}$.
This contradicts Lemma \ref{jisuukeisan},
and hence we obtain $S^1 = S^2 \neq D_j$ for $j=2,3$ with $p_j = 1$.
This implies that
$c_2^*$ and $r_1^*$ are not adjacent and
$c_3^*$ and $r_i^*$ are not adjacent for $i=1,2$.

Let
\[
\begin{array}{l}
\mbox{$\hat{E} = \{a_1, a_3\}$,  $\hat{c}_{1,2} =b_1$  and $\hat{C}_{1} = C_{1,2}$ in the case $R^1 = C_{1,1}$ and $p_1 = 2$,}\\
\mbox{$\hat{E} = \{a_1, a_3\}$,  $\hat{c}_{1,2} =b_1'$ and $\hat{C}_{1} = C_{1,2}$ in the case $R^1 = C_{1,1}$ and $p_1 = 3$,}\\
\mbox{$\hat{E} = \{b_1, b_3\}$,  $\hat{c}_{1,2} =a_1$  and $\hat{C}_{1} = C_{1,1}$ in the case $R^1 = C_{1,2}$ and $p_1 = 2$ and}\\
\mbox{$\hat{E} = \{b_1', b_1\}$, $\hat{c}_{1,2} =a_1$  and $\hat{C}_{1} = C_{1,1}$ in the case $R^1 = C_{1,2}$ and $p_1 = 3$.}\\
\end{array}
\]
In either case,
$\{r_1, s_1, c_{1,1}, \hat{c}_{1,2}\} \cup \hat{E}$ induces a subdivided claw.
By Corollary \ref{endmove} i) and ii),
there exist a heavy edge $c_1^{**} \in \{\hat{c}_{1,2}, c_{1,1}\} \cup l(\hat{C}^1) \cup l(R^1)$.
If $c_1^{**} \in \{\hat{c}_{1,2}\} \cup l(\hat{C}_1)$,
then $\{c_1^{**}, r_1^*, c_2^*, c_3^*\}$ is a heavy matching,
a contradiction.
Hence $c_1^{**} \in \{c_{1,1}\} \cup l(R^1)$.
By Corollary \ref{endmove} ii), $\ed(c_1^{**}) \ge \frac{n-2}3+2$.
Then $\{c_1^{**}, r_2^*, c_3^*\}$ is a matching with
$\ed(c_1^{**}) + \ed(r_2^*) + \ed(c_3^*) \ge n+2$,
which contradicts Lemma \ref{jisuukeisan}.
\qed
\medskip

Without loss of generality, we may assume that
$R^1 \in \{C_{1,1}, C_{1,2}\}$.
Note that
Claim \ref{hachigatu_new} yields
\begin{eqnarray}
\mbox{$V(Q_{u}') \cap V(T) = \{R^1\}$
for any path $Q_{u}'$ which joins an endvertex of $u$ and a vertex in $T$.}
\label{onlyloop}
\end{eqnarray}
Take $r_1^*$ as in Claim \ref{siheavy}.
Since $Q_u$ is a closed trail and $H$ is a triangle-free simple graph,
$|E(Q_u)| \ge 4$.
Hence we can take $r_3 \in E(Q)$ so that $r_3$ and $r_1^*$ are not adjacent.
By (\ref{onlyloop}),
$r_i$ and $c_j^*$ are not adjacent for each $i \in \{1,3\}$ and $j \in \{2,3\}$.
This implies that
$\{r_1^*, c_2^*, c_3^*\}$ is a matching in $H$.

Again by (\ref{onlyloop}),
for each $i\in \{2,3\}$,
neither of the two endvertices of $c_i^*$ is adjacent to $S^1$.
Moreover, since $H$ is triangle-free,
$R^1$ is adjacent to at most one of the endvertices of $c_i^*$.
Hence
\begin{eqnarray}
\mbox{$|N_H^e(r_1^*) \cap N_H^e(c_i^*)| \le 1$ for $i=2,3$.}
\label{aida}
\end{eqnarray}
Let $\gamma = \sum_{\{e_1, e_2\} \subseteq \{r_1^*, c_2^*, c_3^*\}} |N_H^e(e_1) \cap N_H^e(e_2)|$.
Since $r_3 \notin N_H^e(r_1^*) \cup N_H^e(c_2^*) \cup N_H^e(c_3^*)$,
\begin{eqnarray*}
n = |E(H)|
&\ge& |N_H^e(r_1^*)| + |N_H^e(c_2^*)| + |N_H^e(c_3^*)| + |\{r_1^*, c_2^*, c_3^* \}| - \gamma + |\{r_3\}|\\
&\ge& \frac{n-2}3 + 2 +  \frac{n-2}3 + \frac{n-2}3 + 3 - \gamma + 1 \; = n+4 - \gamma,
\end{eqnarray*}
and hence $\gamma \ge 4$.
Since $H$ is triangle-free,
$|N_H^e(e_1) \cap N_H^e(e_2)| \le 2$
for every pair of non-adjacent edges $e_1, e_2 \in E(H)$.
Hence by (\ref{aida}),
we have
$|N_H^e(c_2^*) \cap N_H^e(c_3^*)| = 2$,
$|N_H^e(r_1^*) \cap N_H^e(c_i^*)| = 1$ for $i=2,3$,
and
$E(H) = N_H^e(r_1^*) \cup N_H^e(c_2^*) \cup N_H^e(c_3^*) \cup \{r_3\}$.
This yields $p_1 = 1$,
since otherwise $a_1$ or $b_1$ is not contained in
$N_H^e(r_1^*) \cup N_H^e(c_2^*) \cup N_H^e(c_3^*) \cup \{r_3\}$.

Let $i \in \{2,3\}$.
If $p_i = 3$,
then $b_i$ is not contained in
$N_H^e(r_1^*) \cup N_H^e(c_2^*) \cup N_H^e(c_3^*) \cup \{r_3\}$,
a contradiction.
Hence $p_i \le 2$.
Recall that
neither of the endvertices of $c_i^*$ is adjacent to $S^1$.
Since $|N_H^e(r_1^*) \cap N_H^e(c_i^*)| = 1$,
$R^1$ has a neighbor in an endvertex of $c_i^*$.
On the other hand, $p_1 = 1$ yields $R^1 = C_{1,1}$,
and since $H$ is triangle-free, $C_{1,1}C_{i,1}, C_{1,1}C_{i,2} \notin E(H)$.
Thus we can deduce that $p_i = 1$, $c_i^* = c_{i,1}$ and
$C_{1,1}D_i \in E(H)$ for $i=2,3$.

Let $Q'$ be the path $C_{2,1} D_2 C_{1,1} D_3 C_{3,1}$,
then $T' = P_1 \cup P_2 \cup P_3 \cup Q_u \cup Q' - \{a_2, b_3\}$
is a closed trail passing through $c_2^*$, $c_3^*$ and $r_3$.
Moreover, $T'$ contains $V(r_1)$.
Since $E(H) = N_H^e(r_1^*) \cup N_H^e(c_2^*) \cup N_H^e(c_3^*) \cup \{r_3\}$ and
$r_1^* \in \{r_1\} \cup l(R_1)$,
$T'$ is a DCT of $H$,
a contradiction.
This completes the proof of Case 1.
\medskip\\
%
%
%
%
%
%
%
%
%
%
\noindent
\textbf{Case 2.}\quad
$c_i^* \in \{a_i\} \cup l(A)$ for some $i$ and 
$\bar{c}_j^* \in \{b_j\} \cup l(B)$ for some $j$.
\medskip

Recall that, by Corollary \ref{endmove} i) and ii),
$\ed(c_i^*), \ed(\bar{c}_j^*) \ge \frac{n-2}3+2$.
\medskip
\\
\textbf{Subcase 2.1.}\quad
$c_k^* \in \{c_{k,1}\} \cup l(C_{k,1})$ or
$\bar{c}_k^* \in \{b_k'\} \cup l(C_{k, p_k})$
holds for some $k \in \{1,2,3\}\setminus \{i,j\}$.
\medskip

Without loss of generality, we may assume that 
$c_k^* \in \{c_{k,1}\} \cup l(C_{k,1})$.
If $i \neq j$, then
$\{c_i^*, \bar{c}_j^*, c_k^*\}$ is a heavy matching with
$\ed(c_i^*) + \ed(\bar{c}_j^*) + \ed(c_k^*) \ge
2\left(\frac{n-2}3 + 2\right) + \frac{n-2}3 =n+2$,
which contradicts Lemma \ref{jisuukeisan}.
Hence we have $i = j$.
Without loss of generality,
we may assume that $i=1$ and $k = 3$.
Then we can deduce that
$c_2^* \in \{c_{2,1}\} \cup l(C_{2,1})$,
for otherwise Corollary \ref{endmove} ii) implies
$\ed(c_2^*) \ge \frac{n-2}3 +2$, and so
$\{\bar{c}_1^*, c_2^*, c_3^*\}$ is a heavy matching with
$\ed(\bar{c}_1^*) + \ed(c_2^*) + \ed(c_3^*) \ge n+2$.
Moreover,
we have $p_1 = 1$,
for otherwise 
$\{c_1^*, \bar{c}_1^*, c_2^*\}$ is a heavy matching with
$\ed(c_1^*) + \ed(\bar{c}_1^*) + \ed(c_2^*) \ge n+2$.

Let
$\Xi = H[\bigcup_{1\le i  \le 3} V(c_i^*)]$,
$E_0 = E(H-V(\Xi))$ and
$\Gamma = \bigcup_{1\le i < j \le 3} (N_H^e(c_i^*) \cap N_H^e(c_j^*))$.
Then
$n = |E(H)| \ge |N_H^e(c_1^*)| + |N_H^e(c_2^*)| + |N_H^e(c_3^*)| + |\{c_1^*, c_2^*, c_3^* \}| - |\Gamma| + |E_0|
\ge \frac{n-2}3 + 2 +  \frac{n-2}3 + \frac{n-2}3 + 3 - |\Gamma| + |E_0| \; = n+3 - |\Gamma| + |E_0|$,
thus
\begin{eqnarray}
|\Gamma| \ge |E_0|+3.
\label{eq_case2}
\end{eqnarray}
Since $H$ is triangle-free, we have $|\Gamma| \le 6$,
and hence we have $|E_0| \le 3$.
If $\Xi$ is collapsible,
then we obtain a DCT of $H$ by Corollary \ref{nokorisanhen_c},
a contradiction.
Hence $\Xi$ is not collapsible.
By Lemma \ref{lm_k33-} and the fact that $H$ is triangle-free,
we have $|\Gamma| \le 4$.
Hence it follows from (\ref{eq_case2}) that
$|E_0| \le 1.$

Let
$C_i = D_i$ in the case $p_i = 1$ and $C_i = C_{i,2}$ in the case $p_i \ge 2$.
Since $c_1^* \in \{a_1\} \cup l(A)$
and 
$c_i^* \in \{c_{i,1}\} \cup l(C_{i,1})$ for $i=2,3$,
$\Gamma \subseteq E(H[V(a_1) \cup V(c_{2,1}) \cup V(c_{3,1})])$.
Thus each edge of $\Gamma$ joins two vertices of 
$\{A, C_{1,1}, C_{2,1}, C_2, C_{3,1}, C_3\}$.
On the other hand,
since $H$ is triangle-free,
$A C_2, A C_3, C_{1,1} C_{2,1}, C_{2,1} C_{3,1}, C_{3,1} C_{1,1} \notin E(H)$.
Hence, for every $e \in \Gamma \setminus \{A C_{2,1}, A C_{3,1}\}$, $e = C_i C_{j,1}$ or
$C_i C_j$ with $i \neq 1,j$.
\begin{claim}
$p_i \le 2$ for $i=2, 3$.
\label{Augest24_4}
\end{claim}
\noindent\textit{Proof.}\quad
Assume not.
By symmetry, we may assume that $p_3 \ge 3$.
Then the fact $|E_0| \le 1$ yields
$p_2 \le 2$, $p_3 = 3$ and $E_0 = \{b_3\}$.

Assume that $C_{1,1} C_{3,2} \in E(H)$.
Then $\{C_{1,1} C_{3,2}, c_{1,1}, c_{3,1}, a_3, c_{3,2}, b_3\}$ induces
a subdivided claw, and hence
there exists a heavy edge $b_3^{**} \in \{b_3, c_{3,2}\} \cup l(C_{3,3}) \cup l(C_{3,2})$
by Corollary \ref{endmove} i) and ii).
If $b_3^{**} \in \{b_3\} \cup l(C_{3,3})$,
then $\{c_1^*, c_2^*, c_3^*, b_3^{**}\}$ is a heavy matching of size $4$,
a contradiction.
Moreover, if $b_3^{**} \in \{c_{3,2}\} \cup l(C_{3,2})$,
then since Corollary \ref{endmove} ii) yields $\ed(b_3^{**}) \ge \frac{n-2}3 + 2$,
$\{c_1^*, c_2^*, b_3^{**}\}$ is a heavy matching with
$\ed(c_1^*) + \ed(c_2^*) + \ed(b_3^{**}) \ge n+2$,
which contradicts Lemma \ref{jisuukeisan}.
Hence $C_{1,1} C_{3,2} \notin E(H)$.

If $C_2 C_{i,1} \in E(H)$ for $i=1$ or $3$,
then
$T = A C_{2,1} C_2 C_{i,1} \cup (P_i-a_i) \cup P_j$,
where $j \in \{1,2,3\} \setminus \{i,2\}$,
is a closed trail containing all the vertices in
$\{A, C_{1,1}, C_{2,1}, C_2, C_{3,1}, C_{3,2}, B\}=
V(a_1) \cup V(c_{2,1}) \cup V(c_{3,1}) \cup \{B\}$.
Since $b_3$ is dominated by $B$,
$T$ is a DCT of $H$,
a contradiction.
Hence $C_2 C_{i,1} \notin E(H)$ for $i=1$ and $3$.

By (\ref{eq_case2}), we have $|\Gamma| = 4$.
Since $C_{1,1} C_{3,2}, C_2 C_{1,1}, C_2 C_{3,1} \notin E(H)$,
$\Gamma = \{A C_{2,1}, A C_{3,1},$ $C_{2,1} C_{3,2}, C_2 C_{3,2}\}$.
Then $C_{2,1} C_2 C_{3,2}$ is a triangle, a contradiction.
\qed
\medskip

By Claim \ref{Augest24_4},
we obtain $X=V(P_1 \cup P_2 \cup P_3)$.
\begin{claim}
There exists a closed trail $T$ of $H$
such that $X \subseteq V(T) \subseteq X \cup \{D_i \mid p_i = 1\}$.
\label{Augest24_4_2}
\end{claim}
\noindent\textit{Proof.}\quad
Since (\ref{eq_case2}) yields $|\Gamma| \ge 3$,
there exists an edge $e \in \Gamma \setminus \{A C_{2,1}, A C_{3,1}\}$.
If $C_i C_{j,1} \in E(H)$ for some $i \neq j$, then
$A C_{i,1} C_i C_{j,1} \cup (P_j - a_j) \cup P_k$,
where $k \in \{1,2,3\} \setminus \{i,j\}$,
is a required closed trail.
Hence we may assume that
$e = C_2 C_3$.
Then by Claim \ref{Augest24_4}, either $p_2 =1 $ or $p_3 = 1$ holds since $H$ is triangle-free.
Without loss of generality, we may assume that $p_2 = 1$.
Then
$A C_{3,1} C_3 C_2 C_{2,1} B \cup P_1$
is a required closed trail.
\qed
\medskip

\begin{claim}
If $u \in E(H-X) \setminus \{D_2D_3\}$,
then $u$ is not heavy.
\label{Augest24_4_3}
\end{claim}
\noindent\textit{Proof.}\quad
Assume to the contrary that $u$ is heavy.
Since $u \in E(H-X)$,
$u$ and $c_1^*$ are not adjacent in $H$.
Moreover,
if $u$ and $c_i^*$ are adjacent in $H$ for $i=2$ or $3$,
then $p_i = 1$ and $D_i$ is the common endvertex of $u$ and $c_i^*$.
Since $u \neq D_2D_3$,
either $c_2^*$ or $c_3^*$ is not adjacent to $u$.
If $u$ is adjacent to neither $c_2^*$ nor $c_3^*$,
then $\{c_1^*, c_2^*, c_3^*, u\}$ is a heavy matching of size $4$,
a contradiction.
Hence, without loss of generality,
we may assume that $u$ is adjacent to $c_2^*$ but not to $c_3^*$. 
Then
$p_2 = 1$ and $D_2$ is the common endvertex of $u$ and $c_2^*$.

Since $\{a_1, c_{1,1}, a_2, b_2, a_3, c_{3,1}\}$ induces
a subdivided claw in $H$,
by Corollary \ref{endmove} i) and ii),
there exists a heavy edge $b_2^{**} \in \{a_2, b_2\} \cup l(C_{2,1}) \cup l(A)$.
If $b_2^{**} \in \{b_2\} \cup l(C_{2,1})$,
then $\{c_1^*, b_2^{**}, c_3^*, u\}$ is a heavy matching of size $4$,
and if $b_2^{**} \in \{a_2\} \cup l(A)$,
then $\{\bar{c}_1^*, b_2^{**}, c_3^*, u\}$ is a heavy matching of size $4$,
a contradiction.
\qed
\medskip

It follows from Claim \ref{Augest24_4_2}
that $E(H-X) \neq \emptyset$.
Moreover, in the case $p_2 = p_3 = 1$ and $D_2 D_3 \in E(G)$,
we can deduce that
$E(H-X)\setminus \{D_2 D_3\} \neq \emptyset$,
since otherwise $A C_{2,1} D_2 D_3 C_{3,1} B C_{1,1} A$
is a DCT of $H$.
Since $H$ is connected,
we can take $u \in E(H-X)\setminus \{D_2 D_3\}$
so that an endvertex $S$ of $u$ is adjacent to
a vertex $R \in X$
(possibly $S = D_i$ for some $i$; in this case let $R = C_{i,1}$).
Let $S'$ be the other endvertex of $u$ and let $r = SR$.

We shall prove that
there exist two paths $\Lambda_1$, $\Lambda_2$ of length two
such that
$\{u,r\} \cup E(\Lambda_1) \cup E(\Lambda_2)$ induces a subdivided claw.
If $p_i \ge 2$ for some $i$,
then $R$ is contained in a cycle of length at least $5$ in $P_1 \cup P_2 \cup P_3$.
Since $S, S' \notin X$,
we can find $\Lambda_1$ and $\Lambda_2$ in this cycle.
Hence we consider the case where  $p_i = 1$ for each $i$.
If $R = C_{i,1}$ for some $i$,
then we can find $\Lambda_1$ and $\Lambda_2$ from $P_1 \cup P_2 \cup P_3 - \{a_j, b_k\}$,
where $\{i,j,k\} = \{1,2,3\}$.
If $R = A$ or $B$,
then we can take $j,k$ so that $S' \neq D_j, D_k$,
and then 
$\Lambda_1 = R C_{j,1} D_j$ and
$\Lambda_2 = R C_{k,1} D_k$ are
the desired paths.

Since $\{u,r\} \cup E(\Lambda_1) \cup E(\Lambda_2)$ induces a subdivided claw,
by Corollary \ref{endmove} i) and ii),
there exists a heavy edge $u^* \in \{u,r\} \cup l(S) \cup l(R)$.
By Claim \ref{Augest24_4_3},
we obtain $u^* \in \{r\} \cup l(R)$,
and hence $\ed(u^*) \ge \frac{n-2}3 + 2$.
If $R = A$ (resp.~$B$), then
$\{u^*, \bar{c}_1^*, c_2^*, c_3^* \}$
(resp.~$\{u^*, c_1^*, c_2^*, c_3^* \}$)
is a heavy matching of size $4$, a contradiction.
If $R \in V(P_i) \setminus \{A,B\}$ for $i=2$ or $3$,
then $S \neq D_j$ follows from the choice of $R$, where $j \in \{2,3\} \setminus \{i\}$.
Hence $\{u^*, c_1^*, c_j^*\}$ is a matching with
$\ed(u^*) + \ed(c_1^*) + \ed(c_j^*) \ge n+2$,
which contradicts Lemma \ref{jisuukeisan}.
Therefore we have $R = C_{1,1}$.
Note that,
by the above argument,
we can deduce that
\begin{eqnarray}
\mbox{$V(Q_{u}) \cap X = \{C_{1,1}\}$
for any path $Q_{u}$ which joins an endvertex of $u$ and a vertex in $X$.}
\label{onlyloop2}
\end{eqnarray}
Since $H$ is essentially $2$-edge-connected,
we can take two edge-disjoint paths $Q_{u}^1$ and $Q_{u}^2$
each of which joins $C_{1,1}$ and an endvertex of $u$.
By (\ref{onlyloop2}), we obtain
$V(Q_{u}^i) \cap V(\Xi) = \{C_{1,1}\}$
for $i=1,2$.
Since $H$ is triangle-free simple graph,
we can find two edges in $E(Q_{u}^1) \cup E(Q_{u}^2) \cup \{u\}$
which are not dominated by any vertex in $V(\Xi)$.
Hence $|E_0| \ge 2$, a contradiction.
\medskip
\\
\noindent%
\textbf{Subcase 2.2.}\quad
$c_k^* \in \{a_k\} \cup l(A)$ and
$\bar{c}_k^* \in \{b_k\} \cup l(B)$
holds for any $k \in \{1,2,3\}\setminus \{i,j\}$.

\begin{claim}
\mbox{$c_h^* \in \{a_h\} \cup l(A)$ and
$\bar{c}_h^* \in \{b_h\} \cup l(B)$
holds for any $h \in \{1,2,3\}$.}
\label{zenbukuikomu}
\end{claim}
\noindent\textit{Proof.}\quad
Recall that $c_i^* \in \{a_i\} \cup l(A)$ for some $i$ and 
$\bar{c}_j^* \in \{b_j\} \cup l(B)$ for some $j$
by the assumption of Case 2.
If $i = j$, then the claim follows from the assumption of Subcase 2.2.
Hence we assume $i \neq j$.
By Corollary \ref{endmove} i) and ii),
$\ed(c_k^*), \ed(\bar{c}_k^*) \ge \frac{n-2}3 + 2$ for $k \in \{1,2,3\} \setminus \{i,j\}$.
If
$\bar{c}_i^* \in \{b_i'\} \cup l(C_{i,p_i})$
or
$c_j^* \in \{c_{j,1}\} \cup l(C_{j,1})$,
then in the former case
$\{\bar{c}_i^*, \bar{c}_j^*, c_k^*\}$ is a heavy matching
with $\ed(\bar{c}_i^*) + \ed(\bar{c}_j^*) + \ed(c_k^*) \ge n+2$,
and in the latter case
$\{c_i^*, c_j^*, \bar{c}_k^*\}$ is a heavy matching
with $\ed(c_i^*) + \ed(c_j^*) + \ed(\bar{c}_k^*) \ge n+2$.
This contradicts Lemma \ref{jisuukeisan},
and thus 
$\bar{c}_i^* \in \{b_i\} \cup l(B)$
and $c_j^* \in \{a_j\} \cup l(A)$.
Hence the claim holds.
\qed
\medskip

Now we turn our attention to the graphs $G$ and $\cl(G)$.
For $U \in V(H)$,
let $E_U = \{e \in E(H) \mid e \mbox{ is incident with } U\}$.
Then
$E_U \subset V(G)$ and $E_U$ induces a clique in $\cl(G)$.

Let $I_A = \{c_1^*, c_2^*, c_3^*\}$ and
$I_B=\{\bar{c}_1^*, \bar{c}_2^*, \bar{c}_3^*\}$.
By Lemma \ref{em_g}, Corollary \ref{endmove} i), ii) and Claim \ref{zenbukuikomu},
there exist induced nets $N_A$ and $N_B$ of $G$
such that the vertices in $I_A$ are the endvertices of $N_A$ and
the vertices in $I_B$ are the endvertices of $N_B$.
Hence $d_G(v) \ge \frac{n-2}3$ for each $v \in I_A \cup I_B$.
By Claim \ref{zenbukuikomu},
$I_A \subseteq E_A$ and $I_B \subseteq E_B$.
Note that, since there is no eligible vertex in $\cl(G)$,
\begin{eqnarray}
\mbox{$|N_{\cl(G)}(y) \cap V(Z)|\le 1$
if $Z$ is a clique of $\cl(G)$ and $y \notin V(Z)$.}
\label{closure}
\end{eqnarray}

If there exists $z \in E_A \cap E_B$,
then $z$ is an edge of $H$ joining $A$ and $B$.
Since $c_i^* \in \{a_i\} \cup l(A)$ and
$\bar{c}_i^* \in \{b_i\} \cup l(B)$,
we obtain $v \notin E_A \cap E_B$ for every $v \in I_A \cup I_B$.

Assume that there exists $v \in V(G)$ such that
$|N_G(v) \cap (I_A \cup I_B)| \ge 3$.
Since $G$ is claw-free and both of $I_A$ and $I_B$ are independent sets in $G$,
$|N_G(v) \cap I_A|$, $|N_G(v) \cap I_B| \le 2$.
Hence, without loss of generality,
we may assume that  $|N_G(v) \cap I_A| = 2$.
Then $|N_G(v) \cap I_B| \ge 1$.
Since $E(G) \subseteq E(\cl(G))$,
it follows from (\ref{closure}) that $v \in E_A$.
Again by the claw-freeness of $G$,
there exists $v_A \in N_G(v) \cap I_A$ and $v_B \in N_G(v) \cap I_B$
such that $v_A v_B \in E(G)$.
Then $v_A, v \in N_{cl(G)}(v_B)$.
Since $v_A, v \in E_A$,
(\ref{closure}) yields $v_B \in E_A$,
which contradicts the fact that $v_B \notin E_A \cap E_B$.
Therefore we have $|N_G(v) \cap (I_A \cup I_B)| \le 2$ for each $v \in V(G)$.
Furthermore,
(\ref{closure}) yields $|N_G(v) \cap I_B| \le 1$ for each $v \in I_A$
and $|N_G(v) \cap I_A| \le 1$ for each $v \in I_B$,
and hence $|N_G(v) \cap (I_A \cup I_B)| \le 1$ for each $v \in I_A \cup I_B$.
Therefore
\[
\sum_{v \in I_A \cup I_B} d_G(v)
 = \sum_{v \in V(G)} |N_G(v) \cap (I_A \cup I_B)|
 \le 1 \times 6 + 2 \times (n-6) = 2n - 6, 
\]
which contradicts the fact that
$d_G(v) \ge \frac{n-2}3$ for each $v \in I_A \cup I_B$.
\qed
\begin{thm}
Conjecture \ref{conj} is true for graphs with at least $33$ vertices.
\label{atleast33}
\end{thm}
\noindent
\textit{Proof.}\quad
Let $G$ be a graph of order at least $33$
which satisfies the assumption of Conjecture \ref{conj}
and let $H$ be the triangle-free graph such that
$L(H) = \cl(G)$.
Assume that there exists a heavy matching of size $4$,
say $\{e_1, e_2, e_3, e_4\}$, in $H$.
Since $H$ is triangle-free,
$|N_H^e(e_i) \cap N_H^e(e_j)| \le 2$ for each $i \neq j$.
Then
\begin{eqnarray*}
|E(H)| & \ge & \sum_{i=1}^4 \ed(e_i) + |\{e_1, e_2, e_3, e_4\}|
 - \sum_{1 \le i < j \le 4}|N_H^e(e_i) \cap N_H^e(e_j)| \\
& \ge & 4 \cdot \frac{|E(H)|-2}3 + 4  - 6 \times 2
\: = \: \frac{4|E(H)| -32}3, 
\end{eqnarray*}
which yields $|V(G)| = |E(H)|  \le 32$, a contradiction.
Hence there exists no heavy matching of size $4$ in $H$.
Therefore,
by Theorem \ref{main},
there exists a DCT of $H$,
and by Theorems \ref{Ry} and \ref{HN}, $G$ is hamiltonian.
\qed
\section{Sketch of the proof for the small graphs}
\label{small}
In this section we provide a sketch of the proof of Conjecture \ref{conj}
for graphs of order at most $32$. For a detailed proof,
we refer the readers to \cite{Arxiv}.

Let $H$ be the triangle-free graph such that
$L(H) = \cl(G)$.
By Theorems \ref{Ry} and \ref{HN},
it suffices to prove that
$H$ has a DCT.

\begin{figure}[tb]
\begin{center}
\begin{picture}(120,100)
\put(45,100){\footnotesize$x_1$}
\put(143,100){\footnotesize$x_2$}
%
\put(45,65){\footnotesize$w_1$}
\put(143,65){\footnotesize$w_2$}
%
\put(45,24){\footnotesize$u_1$}
\put(143,24){\footnotesize$u_2$}
\includegraphics[width=5cm]{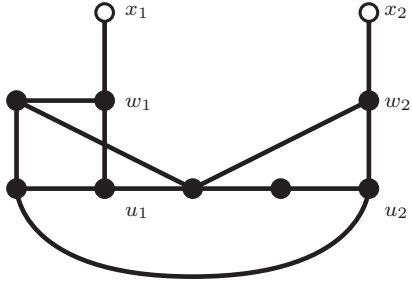}
\end{picture}
\caption{The graph $\Xi^*$ (induced by black vertices).}
\label{fig_Xi}
\end{center}
\end{figure}

First consider the case $n \ge 15$.
By Theorem \ref{main},
we may assume that there exists a heavy matching $M$ of size $4$ in $H$.
Let
$\Xi^* = H[V(M)]$ and
$E_0 = E(H - V(\Xi^*))$.
Since $\Xi^*$ is triangle-free,
we obtain $|E(\Xi^*)| \le 16$.
Moreover, since
$n  =  \sum_{e \in M} \ed(e) + |M| - |E(\Xi^*) \setminus M| + |E_0|
 \ge  4 \cdot \frac{n-2}3 + 4 - (|E(\Xi^*)|-4) +|E_0|$,
we have
\begin{eqnarray}
|E_0| \le |E(\Xi^*)| - \frac{n+16}3.
\label{Dec14}
\end{eqnarray}
Since $H$ is essentially $2$-edge-connected,
we have $|E(H) \setminus \left( E(\Xi^*) \cup E_0\right) | \ge 2$
if $E_0 \neq \emptyset$
(consider two edge-disjoint paths joining an edge of $E_0$ and $\Xi^*$).
This implies
\begin{eqnarray}
|E_0| \le \max\{0, n - |E(\Xi^*)| - 2\}.
\label{Sep29_2}
\end{eqnarray}
If $|E(\Xi^*)| \ge 12$,
then by examining all the possible cases
(note that $\Xi^*$ may not be bipartite in the case $|E(\Xi^*)| = 12$),
we can deduce that either $\Xi^*$ is collapsible or
there exists a vertex $x$ of degree one in $\Xi^*$.
In the former case,
since (\ref{Dec14}) and (\ref{Sep29_2}) yield $|E_0| \le 3$,
we obtain a DCT of $H$ by Corollary \ref{nokorisanhen_c}.
In the latter case,
we have $|E(\Xi^*)| \le 13$, which yield $|E_0| \le 2$.
Moreover, since $|E(\Xi^*)| \ge 12$,
$\Xi^* - \{x\}$ is collapsible.
Then by the similar argument as in Corollary \ref{nokorisanhen_c},
we obtain a DCT of $H$.
Thus we assume $|E(\Xi^*)| \le 11$.
Then (\ref{Dec14}) and the fact $n \ge 15$ yield $|E(\Xi^*)| = 11$ and $E_0 = \emptyset$,
and hence it suffices to find a spanning closed trail of $\Xi^*$.
If $\Xi^*$ is bipartite and $\Xi^*$ has no spanning closed trail,
there must exist a vertex $x$ with $d_{\Xi^*}(x)=1$
and $\Xi' \subseteq \Xi^* - \{x\}$
which is isomorphic to $K_{3,3}$ or $K_{3,3}^-$.
Then we can find a DCT of $H/\Xi'$
by using the fact that $E_0 = \emptyset$,
which yields a DCT of $H$.
If $\Xi^*$ is non-bipartite, then since $H$ is triangle-free and $|E(\Xi^*)| = 11$,
$\Xi^*$ contains an induced cycle $C$ of length $5$.
Let $W = \Xi^* - V(C)$,
then we have $|E(W)| \le 2$ and the number of edges between $W$ and $C$ is $6-|E(W)|$.
By enumerating all the possible structure of $\Xi^*$
and by examining each case carefully,
we can deduce that 
either $\Xi^*$ has a spanning closed trail,
$H$ has a DCT or
$\Xi^*$ is isomorphic to the graph
which is induced by the black vertices in Figure \ref{fig_Xi},
without using Corollary \ref{endmove}.
In the latter case,
since $\Xi^* - \{w_i\}$ has a spanning closed trail
for each $i=1,2$,
we may assume that each $w_i$ has a neighbor $x_i$ in $H - V(\Xi^*)$.
Then, for each $i$, we obtain a subdivided claw containing $x_i w_i u_i$.
By Corollary \ref{endmove} i) and ii),
we obtain an edge joining a white vertex and a black vertex in Figure \ref{fig_Xi},
which yields a DCT of $H$.

Next consider the case $n \le 14$.
Assume that $H$ does not have a DCT,
and take a closed trail $T \subseteq H$ so that
$T$ dominates as many edges as possible.
Then we can take a component $S$ of $H-V(T)$ containing an edge of $E(H-V(T))$.
Let $X=\{x_1, x_2, \ldots ,x_k\}$ be the set of the vertices of $T$
which have a neighbor in $S$,
where $x_1, x_2, \ldots $ appear in this order along $T$, and
let $x_{k+1} = x_1$.
Let $U = V(H) \setminus (V(T) \cup V(S))$,
let $T_i$ be the segment of $T$ between $x_i$ and $x_{i+1}$ and
let $P_i$ be a path of $H$
joining $x_i$ and $x_{i+1}$ whose internal vertices are contained in $S$.
Moreover, let $F_i$ be the set of edges in $H$
joining a vertex in $T_i - \{x_i, x_{i+1}\}$
and a vertex in $T_i \cup U$ 
and let $S_i$ be the set of edges in $H$
which has at least one endvertex in $P_i  - \{x_i, x_{i+1}\}$.
In the case where $T$ is a cycle (that is, each vertex appears exactly once on $T$),
then $F_1, \ldots , F_k, S_i$ are edge-disjoint for each $i$.
Moreover, by the maximality of the number of edges that $T$ dominates,
we have $|F_i| \ge |S_i|$ for each $i$. 
Since $S$ is a non-trivial component, $|S_i| \ge 3$ for each $i$,
and hence
\begin{eqnarray}
14 \ge n \ge \sum_{i=1}^k |F_i| + |S_1|  \ge 3k + 3,
\label{edgecount}
\end{eqnarray}
which yields $k \le 3$.
Note that, in the last inequality,
each of $|F_i|$ and $|S_1|$ is estimated at $3$,
and the set of edges joining $X$ and $U$ is estimated to be empty.
We derive a contradiction by showing that $n \ge 15$.
In both cases $k=2$ and $3$, we can find a subdivided claw
with center $x_i$ containing two edges of $S_i$,
two edges of $F_i$ and two edges of $F_{i+1}$
for each $i$.
By Corollary \ref{endmove} i) and ii)
and close examination of $|F_i|$ and $|S_i|$,
we obtain many edges which is not counted in the last inequality of (\ref{edgecount})
enough to show that $n \ge 15$.
The case where $T$ is not a cycle is basically similar to the above.
By observing the structure of $H$ throughly,
we can find an induced net with center $x_i$ for some $i$.
Then by Corollary \ref{endmove} i) and ii)
we obtain $n \ge 15$.
%
%
%
%
%
%
%
%
%

%
%
%
\end{document}